\documentclass[11pt]{article}
\usepackage[enableskew]{youngtab}
\usepackage{amsfonts,amssymb,amsmath}
\usepackage{stmaryrd}
\usepackage{mathtools}
\usepackage{graphicx}
\usepackage[poly,arrow,curve,matrix]{xy}
\usepackage{wrapfig}
\usepackage{fullpage}
\usepackage{color}

\newcommand{\id}{{\rm id}}

\newcommand{\comp}{{\rm comp}}

\newcommand{\stab}{{\rm stab}}

\newcommand{\Aut}{{\rm Aut}}
\newcommand{\End}{{\rm End}}
\newcommand{\wSimple}{{\sf wSimple}}

\renewcommand{\S}{{\mathfrak S}}

\renewcommand{\k}{\Bbbk}

\newcommand{\N}{{\mathbb N}}
\newcommand{\Z}{{\mathbb Z}}

\newcommand{\A}{{\mathcal A}}
\newcommand{\B}{{\mathcal B}}

\renewcommand{\b}{{\bf b}}

\renewcommand{\1}{{\bf 1}}

\newcommand{\Res}{\underline{\rm Res}}
\newcommand{\Ind}{\underline{\rm Ind}}
\newcommand{\la}{\langle}
\newcommand{\ra}{\rangle}
\newcommand{\ov}{\overline}

\newcommand{\Mod}{{\sf Mod}}
\newcommand{\wMod}{{\sf wMod}}

\newcommand{\Hom}{{\sf Hom}}

\newcommand{\Max}{{\rm Max}}

\newcommand{\tq}{\,|\,}
\newcommand{\proof}{\noindent {\it Proof.$\;$}}

\newcommand{\qed}{\hfill \rule{1.5mm}{1.5mm}}

\input xy
\xyoption{matrix}\xyoption{arrow}
\def\edge{\ar@{-}}
\def\dedge{\ar@{.}}

\newtheorem{theorem}{Theorem}[section]

\newtheorem{remark}[theorem]{Remark}
\newtheorem{example}[theorem]{Example}

\newtheorem{subtheorem}{Theorem}[subsection]
\newtheorem{subproposition}[subtheorem]{Proposition}
\newtheorem{subdefinition}[subtheorem]{Definition}
\newtheorem{sublemma}[subtheorem]{Lemma}
\newtheorem{subexample}[subtheorem]{Example}

\newtheorem{subcorollary}[subtheorem]{Corollary}
\newtheorem{subremark}[subtheorem]{Remark}

\newcommand{\titre}{Weight modules of quantized Weyl algebras}

\begin{document}

\title{{\vspace{-1.5cm} \bf \titre}}
\author{Vyacheslav Futorny, Laurent Rigal, Andrea Solotar}
\date{}
\maketitle

\begin{abstract}
We develop a general framework for studying relative weight representations for certain pairs consisting of an associative algebra and a
commutative subalgebra. Using these tools we describe projective and simple weight modules for quantum Weyl algebras for generic values of deformation parameters. 
We consider two quantum versions: one by Maltsiniotis and the other one by  Akhavizadegan and Jordan. 
\end{abstract}



\tableofcontents


\section*{Introduction and notation}
There is a natural interest in the representation theory of Weyl algebras, both of finite and infinite rank, as the simplest deformations of 
polynomial rings with numerous applications in different areas of mathematics and physics. 
A systematic study of \emph{weight} representations of Weyl algebras was undertaken in \cite{BBF} following the results of 
\cite{Bl}, \cite{Ba}, \cite{BO}, \cite{DGO}, \cite{BB} among the others.
Generalizations of Weyl algebras --infinite rank, generalized, modular, quantum etc.-- were studied in many papers, see for example
\cite{Ba}, \cite{M}, \cite{MT}, \cite{H}, \cite{FGM}, \cite{C}, \cite{FI}, \cite{AJ}, \cite{LMZ}.

Representations of interest in all papers mentioned above are those on which certain commutative subalgebra acts in a locally finite way. 
General
\emph{Harish-Chandra} categories of such representations for a pair $(R, A)$ consisting of an associative algebra $A$ and a commutative 
subalgebra $R$ were
introduced and studied in  \cite{DFO1}. These categories  play a very important role  in representation
theory, in particular in the 
study of Gelfand-Tsetlin modules for the Lie algebra
$gl_n$ \cite{DFO2}, \cite{O}.

In this paper we consider a general
framework for  subcategories of Harish-Chandra categories consisting
of \emph{relative weight} representations for a pair $(R, A)$ on which $R$ is diagonalizable. 

The pairs $(R, A)$ of our interest are so-called of \emph{strongly-fee type}. 
For such pairs we describe projective and simple relative weight representations in Section 1. 
This approach allows to address weight representations of the above mentioned algebras 
in a systematic way. We hope that the technique developed  here will be useful in the study of the representation theory of many types of algebras.  

As an application we focus in this paper  on the representations of quantum Weyl algebras, motivated by their strong 
connections with representations of quantized universal enveloping algebras, --see \cite{LZZ},  \cite{FKZ}-- and quantum affine Lie 
algebras \cite{FHW}. We consider 
two versions of quantum Weyl algebras.
The first one was introduced 
by G. Maltsiniotis (see \cite{M}) in his attempt to develop a quantum differential calculus on the 
standard deformation of the affine space. 
Though its definition makes it a very natural ring of skew differential operators on the quantum 
affine space (where classical partial derivatives are replaced by so-called Jackson operators, which 
are multiplicative analogues) it has the drawback of being a non-simple algebra, hence loosing a 
major property of the classical Weyl algebra. This observation lead to the
seek of a simple localisation of this algebra and was the original motivation for the works 
\cite{AJ} and \cite{J} of M. Akhavizadegan and D. Jordan. In these works, the authors introduced  a 
a slightly different deformation of the Weyl algebra, technically easier to handle and then 
showed that the two quantum analogues actually share a common simple localisation which, in a sense, 
is a better quantum analogue to the classical Weyl algebra.
These algebras depend on an $n$-tuple of  deformation parameters  $\bar{q}=(q_{1}, \ldots, q_n)\in(\k^\ast)^n$  and an $n \times n$ skew symmetric 
matrix 
$\Lambda=(\lambda_{ij})$ for some $n$. 

 Denote by $B_{n}^{\bar{q},\Lambda}$ the common simple localisation of both versions of
quantum Weyl algebras above. For a certain  commutative
Laurent polynomial subalgebra $R^{\circ}$ of $B_{n}^{\bar{q},\Lambda}$ the pair ($R^{\circ}, B_{n}^{\bar{q},\Lambda}$) is of strongly-fee type. 
We apply the above results to describe the projective $B_{n}^{\bar{q},\Lambda}$-modules in Section 3. In Section 4 we obtain one of our main results: a classification of simple relatively weight $B_{n}^{\bar{q},\Lambda}$-modules for generic values of deformation parameters 
and special matrix $\Lambda$ with all entries equal to $1$ (Theorem \ref{thm-simple}). In order to deal with 
the general matrix $\Lambda$ we use  Zhang twists \cite{Z} and show in Section 5 how the general case can be reduced to the special case when all  entries of matrix $\Lambda$ equal   to $1$. It leads in Section 6 to our main
 result - classification and construction of simple weight 
$B_{n}^{\bar{q},\Lambda}$-modules. We note that the case of rank two quantum Weyl algebra and arbitrary deformation parameters was considered in [\cite{H}, Theorem 6.14].

We fix a field $\k$. Whenever the symbol $\otimes$ is used, it stands for the tensor product over $\k$ of $\k$-vector 
spaces, while $M \otimes_A N$ will denote the tensor product of a right module $M$ and a left module $N$ over a ring $A$ and
if $A$ is a $\k$-algebra, $\Aut_\k(A)$ denotes its group of $\k$-algebra automorphisms.

We will denote by $\mathbb{N}$ the set of non negative integers and by $\mathbb{N}^{\ast}$ the subset of positive integers.

\section*{Acknowledgments}
The authors were supported by the MathAMSud grant "Representations, homology and Hopf algebra". 
V.F. was partially supported by the CNPq grant (304467/2017-0) and by the
CNPq grant (200783/2018-1). A.S. was partially supported by PIP-CONICET 11220150100483CO and UBACyT 20020170100613BA. She is a research member
of CONICET (Argentina).

\section{A general setting for weight categories} \label{section-weight-categories}

The aim of this section is to study the following context, which turns out to appear very often in 
representation theory. Consider a $\k$-algebra $A$ containing a commutative subalgebra $R$. One may 
consider (left) $A$-modules on which the subalgebra $R$ acts in a semisimple way. That is, 
$A$-modules $M$ which can be decomposed as a direct sum of subspaces indexed by characters of $R$  
on which any element of $R$ acts by scalar multiplication by its image by the relevant character. 
Clearly, such {\em generalised eigenspaces} must be $R$-submodules and, if the connection between $A$ and $R$ is close enough, such a decomposition can provide a good understanding of the $A$-module
structure of $M$.

We are particularly interested by the case in which $(A,R)$ is a {\em strongly-free} extension. 
Roughly
speaking, this means that there is a subset $\b$ of $A$ which is a basis of $A$ both
as a right and as a left module with the additional property that commutation relations between $R$ 
and elements of $\b$ are controlled by automorphisms of the $\k$-algebra $R$.

In the first subsection, we recall useful results on weight modules over a commutative $\k$-algebra.
These are modules on which $R$ acts in a semisimple way --in the above sense. In the second 
subsection, given a pair $(A,R)$ consisting of a $\k$-algebra $A$ and a commutative subalgebra $R$, 
we introduce the notion of relative weight $A$-module, that is $A$-modules which are weight modules 
relative to $R$. In particular, we show that if the extension $R \subseteq A$ is 
{\em strongly-free}, then the natural induction from $R$ to $A$ preserves projective modules. This
provides a natural process to construct $A$-modules which are projective in the category of relative 
weight $A$-modules as we show in Subsection 3. Subsection 4 then is devoted to the study of the 
simple objects in the above category. It turns out that they all are quotients of the 
aforementioned projective objects. The last subsection is of a more technical nature: we focus on 
tensor products of strongly free extensions and study how the projective and simple objects discussed above behave in this case.

\subsection{Weight modules over commutative algebras}

In this subsection, unless otherwise specified, $R$ denotes a commutative $\k$-algebra. 

\begin{subdefinition} -- 
A {\em weight} of $R$ is a morphism of $\k$-algebras from $R$ to $\k$. The set of weights of $R$ will be denoted $\widehat{R}$.  
\end{subdefinition}

Let $\Max(R)$ be the set of maximal ideals of $R$. There is a map
\[
\begin{array}{rcl}
\widehat{R} & \longrightarrow & \Max(R) \cr 
\phi & \mapsto & \ker(\phi). 
\end{array}
\]
The link between a weight of $R$ and its kernel will be useful in the sequel. We make it clear in the following remark. 

\begin{subremark} -- \rm \label{poids-comme-projecteurs}
Recall the natural morphism of algebras $\k \longrightarrow R$, $\lambda \mapsto \lambda.1_R$.
\begin{enumerate}
\item Notice that, for all $r\in R$ and for all $\phi \in \widehat{R}$, the element $r-\phi(r).1_R$ belongs to $\ker(\phi)$. From this, it follows at once that
\[
R = \k.1_R \oplus \ker(\phi). 
\]
In particular, the map $\widehat{R} \longrightarrow  \Max(R)$ mentioned above is injective. 
Of course, it does not need be surjective. 
\item Let now $\sigma$ be an automorphism of the $\k$-algebra $R$. 
\begin{enumerate}
\item It is straightforward to verify that $\ker(\phi)=\sigma(\ker(\phi\circ\sigma))$.
\item It follows, using the first item, that
\[
\sigma(\ker(\phi)) = \ker(\phi) \Longleftrightarrow \phi\circ\sigma=\phi. 
\]
\end{enumerate}
\end{enumerate}
\end{subremark}

\begin{subremark} -- \rm \label{action-automorphismes-sur-poids}
It is clear that there is a right action of the group $\Aut_\k(R)$ of $\k$-algebra automorphisms of $R$ on $\widehat{R}$
as follows~:
\[
\begin{array}{rcl}
\widehat{R} \times \Aut_\k(R) & \longrightarrow & \widehat{R} \cr
(\phi,\sigma) & \mapsto & \phi\circ\sigma.
\end{array}
\]
\end{subremark}

\begin{subexample} -- \rm
The case where $R$ is a polynomial algebra or a Laurent polynomial algebra will be of particular importance. 
\begin{enumerate}
\item Let $R=\k[z_1,\dots,z_n]$, with $n\in\N^\ast$. There is a 
commutative diagram
\[
\xymatrix@!C{
\widehat{R} \ar@{->}[rr]^{1:1}\ar@{->}[dr]^{} & & \k^n\ar@{->}[dl]^{} \\
 & \Max(R)& \\
}
\]
where the horizontal map sends $\phi\in\widehat{R}$ to the $n$-tuple $(\phi(z_1),\dots,\phi(z_n))$, 
while the left hand side map sends $\phi\in\widehat{R}$ to its kernel and the right hand side map sends 
$(\alpha_1,\dots,\alpha_n)$ to $\langle z_1-\alpha_1,\dots,z_n-\alpha_n \rangle$. \\
Clearly, the horizontal map is bijective, while the other two are injective
(see Remark \ref{poids-comme-projecteurs}).
\item There is a similar diagram for the Laurent polynomial algebra in the indeterminates 
$z_1,\dots,z_n$, $n\in\N^\ast$, replacing $\k$ by $\k^\ast$.
\end{enumerate}
\end{subexample}

\begin{subdefinition} -- 
Let $M$ be an $R$-module. Given $\phi\in\widehat{R}$, we say that $m\in M$ is an {\em element of weight} $\phi$ if, for all
$r\in R$, $r.m = \phi(r)m$. We say that an element $m\in M$ is a {\em weight element} if there exists $\phi\in\widehat{R}$
such that $m$ is an element of weight $\phi$. Further, for $\phi\in\widehat{R}$, we set 
\[
M(\phi) = \{m\in M \tq r.m=\phi(r)m, \, \forall r\in R\}.
\]
\end{subdefinition}

\begin{subremark} -- \rm  
Let $f: M \longrightarrow N$ be a morphism of $R$-modules. 
For all $\phi\in\widehat{R}$, $f(M(\phi)) \subseteq N(\phi)$. 
\end{subremark}

\begin{subproposition} -- \label{somme-directe} Let $M$ be an $R$-module. 
\begin{enumerate}
\item For all  $\phi\in\widehat{R}$, $M(\phi)$ is an $R$-submodule of $M$.
\item Let $n \in\N^\ast$.
If $\phi_1,\dots,\phi_n$ are pairwise distinct 
elements of $\widehat{R}$ and, for $1 \le i \le n$, $m_i\in M(\phi_i)$, then 
$\sum_{1\le i \le n} m_i = 0$ if and only if $m_i=0$, for all $1 \le i \le n$. 
That is, the sum $\sum_{\phi\in\widehat{R}} M(\phi)$ is direct.
\end{enumerate}
\end{subproposition}

\proof The proof of the first item is straightforward. Let $n \in\N^\ast$.
We prove the second statement by induction on $n$, the case $n=1$ being trivial. 
Assume the results holds for some $n\in\N^\ast$. 
Suppose now that $\phi_1,\dots,\phi_{n+1}$ are pairwise distinct 
elements of $\widehat{R}$ and, for $1 \le i \le n+1$, consider elements $m_i\in M(\phi_i)$ such that 
$\sum_{1\le i \le {n+1}} m_i = 0$. For all $r \in R$, we have that 
\[0 = r.\left(\sum_{1\le i \le {n+1}} m_i\right) = \sum_{1\le i \le {n+1}} r.m_i = \sum_{1\le i \le {n+1}} \phi_i(r).m_i.
\]
It follows that, for all $r \in R$, 
\[
\sum_{1\le i \le n} (\phi_{n+1}(r)-\phi_i(r)).m_i = 0.
\]
The inductive hypothesis yields that, for all $r\in R$ and $1 \le i \le n$, $(\phi_{n+1}(r)-\phi_i(r)).m_i=0$.
But, the $\phi_i$'s being pairwise distinct, we have $m_1= \dots = m_n=0$ and therefore $m_{n+1}=0$.
This finishes the induction. The result is thus proved.\qed

\begin{subdefinition} -- 
Let $M$ be an $R$-module. For $\phi\in\widehat{R}$, we will call $M(\phi)$ the {\em weight subspace} of $M$ of weight $\phi$. Further, 
we say that $M$ is a {\em weight module} provided it is the (direct) sum of its weight subspaces  
\[
M =  \bigoplus_{\phi\in\widehat{R}} M(\phi)
\]
We denote by $R-\wMod$ the full subcategory of $R - \Mod$ whose objects are the weight $R$-modules.
\end{subdefinition}

\begin{subproposition} -- \label{sous-objets-et-poids}
Let $M$ be an $R$-module and let $L$ be an $R$-submodule of $M$. 
\begin{enumerate}
\item For all $\phi\in\widehat{R}$, $L(\phi)=L \cap M(\phi)$. 
\item  If $M$ is a weight module, then so is $L$.
\end{enumerate}
\end{subproposition}

\proof The first item is obvious. Next, suppose $m$ is a non zero element of $L$. There exists
$n\in\N^\ast$ and pairwise distinct elements $\phi_1,\dots,\phi_n$ of $\widehat{R}$ such that 
\begin{equation}\label{m}
m = \sum_{1\le i \le n} m_i 
\end{equation}
with $m_i\in M(\phi_i)\setminus\{0\}$ for all $1 \le i \le n$. So,
\begin{equation}\label{truc}
\forall r\in R, \qquad r.m - \phi_n(r) m = \sum_{1 \le i \le n} (\phi_i(r)-\phi_n(r)) m_i \in L. 
\end{equation}
Further, if $n > 1$, then for all $1 \le i \le n-1$, we may choose $r_i \in R$ such that 
$\phi_i(r_i) \neq \phi_n(r_i)$, and hence 
$r_i.m - \phi_n(r_i) m \in L\setminus\{0\}$ by Proposition \ref{somme-directe}. 

We are now ready to prove the second item. 
Suppose $L$ is not a weight submodule. This means that there exists  
a non zero element in $L$ such that at least one of its weight summands is not in $L$. 
Among such elements, let us choose $m$ with a minimal number $n$
of non zero weight summands. Write $m$ as in (\ref{m}) in such a way that $m_n \notin L$. 
Notice in addition that we must have $n > 1$.
Due to the minimality of $n$, the argument above based on (\ref{truc}) shows that, for all 
$1 \le i \le n-1$, 
$m_i \in L$. It follows that $m_n \in L$. This is a contradiction. Hence $L$ must be a weight 
submodule.\qed

\medskip

Next, we consider the case where $R$ is a tensor product of commutative algebras, that is
$R= R_1\otimes \dots \otimes R_s$ for some $s\in \N^\ast$, where $R_1,\dots,R_s$ are 
commutative $\k$-algebras.

\begin{subproposition} \label{weight-in-tpa}
Let $s\in\N^\ast$, let $R_1,\dots,R_s$ be commutative 
$\k$-algebras and consider the $\k$-algebra $R = R_1\otimes\dots\otimes R_s$.    
\begin{enumerate}
\item For $1 \le i \le s$, let $\phi_i\in\widehat{R_i}$. The $\k$-linear map
\[
\begin{array}{ccrcl}
\phi_1 \dots \phi_s & : & R & \longrightarrow & \k \cr
 & & x_1 \otimes\dots\otimes x_s & \mapsto & \phi_1(x_1) \dots \phi_s(x_s)
\end{array}
\]
is a morphism of $\k$-algebras. Furthermore, we have that
\[
\ker(\phi_1 \dots \phi_s) = W 
\]
where 
\[
W=\ker(\phi_1) \otimes R_2 \otimes\dots\otimes R_s 
+ R_1 \otimes \ker(\phi_2) \otimes\dots\otimes R_s + \dots + R_1 \otimes\dots\otimes R_{s-1} \otimes \ker(\phi_s).
\]
\item The map 
\[
\begin{array}{ccrcl}
 & & \widehat{R_1} \times\dots\times \widehat{R_s} & \longrightarrow & \widehat{R} \cr
 & & (\phi_1 , \dots , \phi_s) & \mapsto & \phi_1 \dots \phi_s
\end{array}
\]
is bijective.
\end{enumerate}
\end{subproposition}

\proof
In order to prove the first statement, set $\phi=\phi_1 \dots \phi_s$. It is easy to see that $W \subseteq \ker(\phi)$.  
According to Remark \ref{poids-comme-projecteurs}, given $1 \le i \le s$, $R_i = \k.1_{R_i} \oplus \ker(\phi_i)$. 
Hence, $R = \k.1_R \oplus W$. It follows that if $r$ is an element of $R$, then there exist $\lambda\in\k$ and $w \in W$ 
such that $r= \lambda.1_R + w$, and, since $W \subseteq \ker(\phi)$, we have $\lambda=\phi(r)$. It is now clear that 
if $r\in\ker(\phi)$, then $r \in W$.\\
The previous point gives rise to a map
\[
\begin{array}{ccrcl}
 & & \widehat{R_1} \times\dots\times \widehat{R_s} & \longrightarrow & \widehat{R} \cr
 & & (\phi_1 , \dots , \phi_s) & \mapsto & \phi_1 \dots \phi_s
\end{array} .
\]
Now, for $1 \le i \le s$, identify $R_i$ with its canonical image in $R$ under the injective algebra morphism $r \mapsto
1 \otimes\dots\otimes 1\otimes r \otimes 1 \otimes\dots\otimes 1$ ($r$ in $ith$ position). 
We may define the map
\[
\begin{array}{ccrcl}
 & & \widehat{R} & \longrightarrow & \widehat{R_1} \times\dots\times \widehat{R_s} \cr
 & & \phi & \mapsto & (\phi_{|R_1}, \dots ,\phi_{|R_s}).
\end{array}
\]
The two precedent maps are inverse to each other, so that they are bijective.
\qed

\subsection{Relative weight modules} \label{sous-section-rwm}

In this subsection, we fix a pair $(R,A)$ where $A$ is a $\k$-algebra and $R$ a commutative $\k$-subalgebra of $A$. 

Any $A$-module $M$ is an $R$-module by restriction of scalars. Hence, 
if $\phi\in\widehat{R}$, we may speak of elements of $M$ of weight $\phi$. 

Let $A-\Mod$ be the category of left $A$-modules. We denote by $A-\wMod$ the full subcategory of $A-\Mod$ whose objects are 
weight $R$-modules. We call such modules \emph{relative weight} modules.
We have the following restriction of scalars functors

\[
\xymatrix@!C{
A-\Mod \ar@{->}[rrr]^{\Res} &&& R-\Mod\\
& & & \\
A-\wMod \ar@{->}[uu]^{can.} \ar@{->}[rrr]^{\Res} &&& R-\wMod\ar@{->}[uu]^{can.}\\
}
\]
Of course, there is also an induction functor $\Ind = A \otimes_R -$ from $R-\Mod$ to $A-\Mod$. In order to be able to 
induce weight modules from $R$ to $A$, we need the pair $(R,A)$ to satisfy some additional hypotheses.

\begin{subdefinition} -- \label{definition-strongly-free}
We say that the pair $(R,A)$ is of {\em strongly-free type} if there exists a subset $\b$ of $A$ 
such that:
\begin{enumerate}
\item $1\in \b$;
\item $\b$ is a basis of $A$  as a left $R$-module;
\item $\forall b \in \b$, there exists an algebra automorphism $\sigma_b$ of $R$ such that, $\forall r\in R$, $rb=b\sigma_b(r)$. 
\end{enumerate}
\end{subdefinition}

\begin{subremark} -- \rm
If the pair $(R,A)$ is of strongly-free type, then, there exists a subset of $A$ which is an $R$-basis of $A$ both
as left and as right module over $R$.
\end{subremark}

\begin{sublemma} -- \label{induced-of-weight}
Assume the pair $(R,A)$ is of strongly free type. Let $\phi\in\widehat{R}$ and suppose $M$ is an $R$-module. Then, 
\begin{enumerate}
\item for all $m\in M(\phi)$, for all $b\in\b$, $b \otimes_R m \in (\Ind(M))(\phi\circ\sigma_b)$;
\item if $M$ is a weight $R$-module, then $\Ind(M)$ is an object of $A-\wMod$.
\end{enumerate}
\end{sublemma}

\proof The proof of the first item is straightforward. The second one follows from the first since $\b$ generates $A$ as left
$R$-module.\qed\\

As a consequence of Lemma \ref{induced-of-weight}, 
there is a commutative diagram for the induction functors:
\[
\xymatrix@!C{
A-\Mod  &&& \ar@{->}[lll]^{\Ind} R-\Mod\\
& & & \\
A-\wMod \ar@{->}[uu]^{can.} &&& R-\wMod \ar@{->}[lll]^{\Ind} \ar@{->}[uu]^{can.}\\
}
\]

\begin{subcorollary} -- \label{paire-adjointe}
Assume the pair $(R,A)$ is of strongly free type.
The pair $(\Ind,\Res)$ is an adjoint pair of functors between $A-\wMod$ and $R-\wMod$. In particular, if $M$ is 
a relative weight 
$A$-module and $N$ is a weight $R$-module, then there is an isomorphism of abelian groups as follows:
\[
\Hom_{A-\wMod}(\Ind(N), M) \cong \Hom_{R-\wMod}(N,\Res(M)). 
\]
\end{subcorollary}

\proof It is well known that the pair $(\Ind,\Res)$ is an adjoint pair of functors between $A-\Mod$ and 
$R-\Mod$. The proof follows immediately from this fact. \qed

\begin{subcorollary} -- \label{ind-et-projectifs}
Assume the pair $(R,A)$ is of strongly free type.
The functor $\Ind \, : \, R-\wMod \longrightarrow A-\wMod$ preserves projective objects.
\end{subcorollary}

\proof This is an immediate consequence of Corollary \ref{paire-adjointe} since $\Res$ is  exact.\qed

\subsection{A class of projective weight $A$-modules}

Let $R$ be a commutative $\k$-algebra. Consider $\phi\in\widehat{R}$. 
The $R$-module $R/\ker(\phi)$ is clearly an object of $R-\wMod$. More precisely, we have that
any element of this module is of weight $\phi$, that is  
$R/\ker(\phi)=\left(R/\ker(\phi)\right)(\phi)$. In the sequel, we will denote
it by $\k_\phi$.

\begin{sublemma} -- \label{kn-projectif}
Given $\phi\in\widehat{R}$, the module $\k_\phi$ is a projective object in the category 
$R-\wMod$. 
\end{sublemma}

\proof Consider the diagram
\[
\xymatrix@!C{
   & \k_\phi \ar@{->}[d]^{f} & \\
N_1 \ar@{->}[r]^{\pi} & N_2 \ar@{->}[r] & 0\\
}
\]
in $R-\wMod$ with exact horizontal row. Using Prop. \ref{somme-directe}, 
we get that there is an element $z$ of weight $\phi$ in $N_1$ such that $\pi(z)=f(\ov{1_R})$.
Since $\k_\phi$ is a one dimensional $\k$-vector space, 
there is a unique $\k$-linear map $g \, : \, \k_\phi \longrightarrow N_1$ sending $\ov{1_R}$ to $z$. Thus, the diagram
\[
\xymatrix@!C{
   & \k_\phi \ar@{->}[d]^{f}\ar@{->}[dl]^{g} & \\
N_1 \ar@{->}[r]^{\pi} & N_2 \ar@{->}[r] & 0\\
}
\]
is commutative. It is easy to see that $g$ is actually a map in $R-\wMod$, since $z$ has weight 
$\phi$. (Here, we use Remark \ref{poids-comme-projecteurs}, point 1.) 
Hence, $\k_\phi$ is projective. \qed\\

For the rest of this section, let $(R,A)$ be a pair as in the introduction to section 
\ref{sous-section-rwm}, 
assume $(R,A)$ is of strongly-free type and let $\b$ be a basis of the $R$-module 
$A$ as in Definition \ref{definition-strongly-free}.\\

Let $\phi\in\widehat{R}$. Set $P_\phi = \Ind(\k_\phi)$. We already know that $P_\phi$ is an object of $A-\wMod$.

\begin{subcorollary} -- 
For any $\phi\in\widehat{R}$, the module $P_\phi$ is a projective object in the category $A-\wMod$. 
\end{subcorollary}

\proof This is immediate by Lemma \ref{kn-projectif} and Corollary \ref{ind-et-projectifs}. \qed\\

Now we collect some properties of the $A$-modules $P_\phi$, $\phi\in\widehat{R}$, that will be useful later. 

\begin{subremark} -- \rm\label{structure-P} {\bf Structure of $P_\phi$, $\phi\in\widehat{R}$}\\
Recall that we assume $(R,A)$ is of strongly-free type and let $\b$ be a basis of the $R$-module $A$ as in Definition \ref{definition-strongly-free}. Fix $\phi\in\widehat{R}$. 
For all $r\in R$, we put $\bar{r} = r + \ker(\phi) \in\k_\phi$. 

\begin{enumerate}
\item Fix $b\in\b$. We have that $Rb=bR \subseteq A$. Further, $Rb$ is a direct summand of $A$ 
both as a left and as a right $R$-submodule. Now, consider the following subset of $P_\phi$:
\[
b \otimes_R \k_\phi :=\{b\otimes_R \bar{r}, \, r\in R \} \subseteq P_\phi.
\]
Clearly, $b \otimes_R \k_\phi$ is the left $R$-submodule of $P_\phi$ generated by 
$b \otimes_R \bar{1}$ as well as the $\k$-subspace of $P_\phi$ generated by 
$b \otimes_R \bar{1}$ (see Remark \ref{poids-comme-projecteurs}). On the other hand, it is the image of the natural left $R$-linear map
\[
bR \otimes_R \k_\phi \stackrel{\subseteq \otimes_R \id}{\longrightarrow} P_\phi .
\]
But, $bR$ being a direct summand of $A$ as a right $R$-submodule, the above map must be injective. Hence, it identifies $bR \otimes_R \k_\phi$ and $b \otimes_R \k_\phi$ as left
$R$-modules.

In addition, there are obvious $\k$-linear maps as follows:
\begin{equation}\label{2-applications}
\begin{array}{rcl}
\k_\phi & \longrightarrow & bR \otimes_R \k_\phi \cr
\bar{s} & \mapsto & b \otimes_R \bar{s}
\end{array},
\quad
\begin{array}{rcl}
bR \otimes_R \k_\phi & \longrightarrow & \k_\phi \cr
br \otimes_R \bar{s} & \mapsto & r\bar{s}
\end{array}
\end{equation}
which are inverse to each other. As a consequence, $bR \otimes_R \k_\phi$ is a one
dimensional $\k$-vector space with basis $b \otimes_R \bar{1}$.

At this stage, it is interesting to notice that the maps (\ref{2-applications}) are not left $R$-linear. However, if we denote by $^{\sigma_b}\k_\phi$ the left $R$ module obtained by
twisting the left action of $R$ on $\k_\phi$ by the automorphism $\sigma_b$, then 
the maps (\ref{2-applications}) induce a left $R$-linear isomorphism
\[
bR \otimes_R \k_\phi \longrightarrow \, ^{\sigma_b}\k_\phi.
\]
\item We are now in position to give a very explicit description of $P_\phi$. Indeed, by 
standard results, we have an isomorphism of left $R$-modules as follows: 
\[
P_\phi 
= A \otimes_R \k_\phi 
= \left( \bigoplus_{b\in\b} bR \right) \otimes_R \k_\phi
\cong \bigoplus_{b\in\b} bR \otimes_R \k_\phi. 
\]
Thus, taking into account the first point, we end up with a $\k$-vector space isomorphism
\[
P_\phi 
\cong \bigoplus_{b\in\b} \k . b \otimes_R \bar{1}, 
\]
where, for all $b\in\b$, $\dim_\k (\k . b \otimes_R \bar{1})=1$. In particular, 
$P_\phi$ is nonzero and $\{b \otimes_R \bar{1}, b\in \b\}$ is a basis of $P_\phi$ as
a $\k$-vector space.
\item The weight structure of $P_\phi$, however, is not yet completely clear. 
Indeed, we have that,  
\[
\forall b\in\b, \quad \k.b \otimes_R \bar{1} \subseteq (P_\phi)(\phi\circ\sigma_b).
\]
But, the above inclusion may be strict. Actually, it is the case if and only if there exists
$b'\in\b$, $b\neq b'$, such that $\phi\circ\sigma_b=\phi\circ\sigma_{b'}$. 
We will come back to this problem in Remark \ref{generic-P}
\end{enumerate}
\end{subremark}

\subsection{Simple objects in $A - \wMod$}

All along this section, we fix a pair $(R,A)$ as in the introduction to section
\ref{sous-section-rwm}, 
assume $(R,A)$ is of strongly-free type and let $\b$ be a basis of the $R$-module 
$A$ as in Definition \ref{definition-strongly-free}. Our aim is to study 
the simple objects of $A - \wMod$. 

\begin{subdefinition} -- 
An object $S$ in $A-\wMod$ is {\em simple} provided it is non zero and it has no nontrivial weight $A$-submodule. 
\end{subdefinition}

\begin{subremark} -- \rm By Proposition \ref{sous-objets-et-poids}, 
an object $M$ in $A - \wMod$ is simple if and only if it is simple as an object of $A - \Mod$.
\end{subremark}
 
The following result shows that any simple object in $A - \wMod$ arises as a simple quotient of some 
$P_\phi$, with $\phi\in\widehat{R}$.

\begin{subproposition} -- \label{ubiquite-des-P} \rm 
Given a simple object $S$  of $A-\wMod$, there exists $\phi\in\widehat{R}$ such that $S$ is isomorphic in $A-\wMod$ to a simple quotient of $P_\phi$.
\end{subproposition}

\proof Let $S$ be a simple object of $A-\wMod$. Since $S$ is not zero, there exists 
$\phi\in\widehat{R}$ such that $S(\phi)\neq 0$. 
Let $x$ be a non zero element in $S(\phi)$. There exists a non zero $R$-linear map $R \longrightarrow S$ such that
$1 \mapsto x \in S(\phi)$. But, since $x \in S(\phi)$, the above map induces a nonzero $R$-linear map 
$\k_\phi \longrightarrow S$. 
This shows that $\Hom_{R-\wMod}(\k_\phi,S) \neq 0$. It follows from Corollary \ref{paire-adjointe} that 
$\Hom_{A-\wMod}(P_\phi,S) \neq 0$. Hence, there exists a map $\varphi \, : \, P_\phi \longrightarrow S$ in $A-\wMod$ 
which is non zero. But, of course, the image of $\varphi$ is a non zero subobject of $S$, so $\varphi$ must be surjective.
The rest is clear. \qed\\

Proposition \ref{ubiquite-des-P} shows the importance of understanding the nonzero simple quotients of the modules 
$P_\phi$. The following remark deals with the easiest case.

\begin{subremark} -- \rm \label{generic-P}
Recall the hypotheses of the beginning of this section. Notice that, clearly, $\sigma_1=\id$.
Fix $\phi\in\widehat{R}$. Assume, further, that the map $\b \longrightarrow \widehat{R}$, given by $b \mapsto \phi\circ\sigma_b$ is injective (see Remark \ref{structure-P}). 
It follows from Remark \ref{structure-P} that the weight subspace of $P_\phi$ of weight $\phi$ is:
\[
(P_\phi)(\phi) = \k. 1 \otimes_R (1+\ker(\phi)). 
\]
Let $M$ be a strict weight $A$-submodule of $P_\phi$. Since $P_\phi$ is generated by its weight  
subspace $(P_\phi)(\phi)$ of weight $\phi$, and since the latter is a one dimensional 
$\k$-vector space, we must have $(P_\phi)(\phi) \cap M = \{0\}$, which leads to: 
\[
M \subseteq \bigoplus_{\psi\neq\phi} (P_\phi)(\psi). 
\]
It follows that the sum of all strict submodules of $P_\phi$ in the category $A-\wMod$ is again a strict submodule. 
That is,
there is a strict submodule of $P_\phi$ in $A-\wMod$ maximum among strict submodules of $P_\phi$ in $A-\wMod$. 
We denote this
maximum submodule by $N_\phi$. Of course, the corresponding quotient
\[
S_\phi := P_\phi/N_\phi  
\]
is a simple object in $A-\wMod$ and it is the unique simple quotient of $P_\phi$ in $A-\wMod$.
Note that the above applies to $P_\phi$ seen as an object of $A-\Mod$, by Proposition
\ref{sous-objets-et-poids}. 
\end{subremark}

\subsection{Tensor product of strongly-free extensions} \label{TP-of-SFE}

Let $n\in\N^\ast$. For each $i$, $1 \le i \le n$, let $(R_i,A_i)$ be a pair of strongly-free type and let $\b_i$ be a subset of
$A_i$ as in Definition \ref{definition-strongly-free}. Let 
\[ R = \bigotimes_{1 \le i \le n} R_i, \qquad
A = \bigotimes_{1 \le i \le n} A_i
\]
and identify $R$ canonically with its image in $A$. (In the sequel, we will make extensive use 
of Proposition \ref{weight-in-tpa}.)

Suppose in addition that, for $1 \le i \le n$, $M_i$ is an object of $A_i-\wMod$. 
Clearly, $\bigotimes_{1 \le i \le n} M_i$ is an object of $A-\wMod$.   

\begin{sublemma} -- 
The pair $(R,A)$ is of strongly-free type.
\end{sublemma}

\proof It is clear that the subset $\b = \bigotimes_{1 \le i \le n} \b_i$ of $A$ satisfies the conditions of 
Definition \ref{definition-strongly-free}. \qed

\begin{subremark} -- \rm \label{corps-residuel-et-pt}
Retain the notation above and, for all $1 \le i \le n$, let $\phi_i$ be an element of $\widehat{R_i}$. Set $\phi=\phi_1 \dots \phi_n$. Recall the description of $\ker(\phi)$ given in Proposition \ref{weight-in-tpa}.
There is an obvious surjective morphism  
\[
R \longrightarrow \bigotimes_{1 \le i \le n} \k_{\phi_i} 
\]
in $R-\Mod$ such that its kernel contains $\ker(\phi)$. Hence, there is a surjective morphism in $R-\Mod$~: 
\[ \k_\phi \longrightarrow \bigotimes_{1 \le i \le n} \k_{\phi_i}
\]
which must be an isomorphism since 
its source and target are both one dimensional vector spaces over $\k$. 
\end{subremark}

The following theorem provides a description of the projective objects in $A-\wMod$ as tensor products of projectives.

\begin{subtheorem} -- \label{iso-P-tenseur-P}
Using the previous notation, let $\phi_i$ be an element of $\widehat{R_i}$ for all $1 \le i \le n$, and put 
$\phi=\phi_1 \dots \phi_n$. 
There is an isomorphism 
\[
P_\phi \cong \bigotimes_{1\le i \le n} P_{\phi_i} 
\]
in $A-\wMod$.
\end{subtheorem}

\proof It is enough to show that there is an isomorphism in $A-\Mod$:
\[
P_\phi \cong \bigotimes_{1\le i \le n} P_{\phi_i}.
\]
This is a standard fact; indeed, by classical results, we have an isomorphism of $\k$-vector spaces
\[
\begin{array}{rcl}
\bigotimes_{1 \le i \le n} P_{\phi_i} 
& =     & \bigotimes_{1 \le i \le n} \left( A_i \otimes_{R_i} \k_{\phi_i}\right) \\
& \cong & \left(\bigotimes_{1 \le i \le n} A_i \right) \otimes_R 
\left( \bigotimes_{1 \le i \le n}\k_{\phi_i}\right) \\
& \cong & \left(\bigotimes_{1 \le i \le n} A_i \right) \otimes_R  \k_\phi \\
& \cong & A \otimes_R  \k_\phi \\
& = & P_\phi
\end{array}
\]
where the isomorphism of the second row is given by  
\[ \otimes_{1 \le i \le n} \left( a_i \otimes_{R_i} k_i\right) \mapsto
\left(\otimes_{1 \le i \le n} a_i \right) \otimes_R 
\left(\otimes_{1 \le i \le n} k_i \right)
\]
and the isomorphism of the third row is given by Remark \ref{corps-residuel-et-pt}.\qed\\

Now we analyse simple modules in $A-\wMod$. 

\begin{sublemma} -- \rm \label{shurian-simple}
If, for $1 \le i \le n$, 
$S_i$ is a simple object in $A_i -\wMod$ such that $\dim_\k (S_i)(\phi)\le 1$ for any
$\phi \in \widehat{R_i}$, then $S_1\otimes \ldots \otimes S_n$ is a 
simple object in $A -\wMod$. 
\end{sublemma}

\proof  
Suppose that $\dim_\k (S_i)(\phi)\le 1$ for any  $\phi \in \widehat{R_i}$, $i=1, \ldots, n$. 
It follows from Proposition \ref{weight-in-tpa} that 
the weight spaces of
the $A$-module $S_{1} \otimes\dots\otimes S_{n}$ are also at most one dimensional $\k$-vector 
spaces, generated by elementary tensors.
Now, let $W$ be a nonzero (weight) submodule of $S_{1} \otimes\dots\otimes S_{n}$. 
Then, $W$ contains a nonzero weight element which, by the above, must be an elementary tensor. 
Let $w = w_1 \otimes\dots\otimes w_n$ be such an element. For all $i$ such that 
$1 \le i \le n$, $w_i$ is a non zero element
of $S_{i}$ which, hence, generates the $A_i$-module $S_{i}$. 
Therefore, $w$ generates the $A$-module $S_{1} \otimes\dots\otimes S_{n}$ and  
$W=S_{1} \otimes\dots\otimes S_{n}$. The statement is proved. \qed

\begin{subremark} -- \rm 
If we suppose that ${\rm End}_{A_i}(S)=\k$ for any simple object $S$ in $A_i -\wMod$, 
$i=1, \ldots, n-1$, then 
for any collection of simple objects $S_i$ in $A_i -\wMod$, $i=1, \ldots, n$,
the $A$-module $S_1\otimes \ldots \otimes S_n$ is simple.
Indeed, under this hypothesis, we have that $A_i$ is Schurian and tensor-simple  
\cite{B}, which implies the statement.
\end{subremark}

Let $S$ be a simple object in $A-\wMod$. By Proposition \ref{ubiquite-des-P} 
there exists $\phi\in\widehat{R}$ such that $S$ is a quotient of $P_{\phi}$, that is $S\simeq P_{\phi}/N$. We have the following result.

\begin{subcorollary} -- \label{iso-L-tenseur-L}
Let $\phi\in\widehat{R}$ be such that all weight subspaces of $P_{\phi}$ have dimension 
at most $1$. For $i=1,\dots,n$, let $\phi_i \in \widehat{R_i}$ be such that 
$\phi=\phi_1 \dots \phi_n$.
\begin{enumerate}
\item The map $\b \to \widehat{R}$ given by $b \mapsto \phi\circ \sigma_b$ is injective.
Further, the module $P_{\phi}$ has a unique maximal strict submodule $N_{\phi}$ and unique simple quotient $S_\phi = P_{\phi}/N_{\phi}$.
\item For all $1 \le i \le n$, The map $\b_i \to \widehat{R_i}$ given by $b \mapsto \phi_i\circ \sigma_b$ is injective.
Further, the module $P_{\phi_i}$ has a unique maximal strict submodule $N_{\phi_i}$ and unique simple quotient $S_{\phi_i} = P_{\phi_i}/N_{\phi_i}$.
\item There is an isomorphism 
\[
S_\phi \cong S_{\phi_1} \otimes \dots \otimes S_{\phi_n} 
\]
in $A - \wMod$.
\end{enumerate}

\end{subcorollary}

\proof
 We start by recalling that such $\phi_i$'s exist due to Proposition \ref{weight-in-tpa}.
\begin{enumerate}
\item This follows immediately from Remarks \ref{structure-P} and \ref{generic-P}.
\item We know from Theorem \ref{iso-P-tenseur-P} that there is an isomorphism $\iota: P_{\phi}\to P_{\phi_1}\otimes \dots \otimes P_{\phi_n}$. Hence, since all weight spaces of $P_{\phi}$ have dimension $1$, by Proposition \ref{weight-in-tpa}, the same must hold for 
$P_{\phi_i}$, for $i=1,\dots, n$. The statement follows as in the first point.
\item Observe that for all $i$, $S_{\phi_i}$ must also have weight spaces of dimension at most one. 
Thus, we are in the hypotheses of the previous lemma; so that  
$S_{\phi_1}\otimes \dots \otimes S_{\phi_n}$ is a simple object in $A-\wMod$.

Now, consider the following submodule of $P_{\phi_1}\otimes \dots \otimes P_{\phi_n}$: 
\[
X 
= N_{\phi_1} \otimes P_{\phi_2} \otimes\dots\otimes P_{\phi_n} 
+ P_{\phi_1} \otimes N_{\phi_2} \otimes\dots\otimes P_{\phi_n}
+ \dots 
+ P_{\phi_1} \otimes\dots\otimes P_{\phi_{n-1}} \otimes N_{\phi_n}.
\]
We have an obvious surjective map 
\[
P_{\phi}\simeq P_{\phi_1}\otimes \dots \otimes P_{\phi_n} \twoheadrightarrow S_{\phi_1}\otimes \dots \otimes S_{\phi_n}
\]
whose kernel contains $\iota^{-1}(X)$, hence a surjective map
\[
P_{\phi}/\iota^{-1}(X) \twoheadrightarrow S_{\phi_1}\otimes \dots \otimes S_{\phi_n}.
\]
We now proceed to show that the latter must be an isomorphism.
Indeed, taking Proposition \ref{weight-in-tpa} into account, the hypotheses on the dimensions of 
weight spaces prove that the weights of 
$S_{\phi_1}\otimes \dots \otimes S_{\phi_n}$ are the products $\psi_1\dots \psi_n$, where $\psi_i$ 
is a weight of $P_{\phi_i}$ that is not a weight of 
$N_{\phi_i}$. Analogously, the weights of $X$ are the products $\psi_1\dots \psi_n$, where $\psi_i$ 
is a weight of $P_{\phi_i}$ and, for some
$j\in \{1, \dots, n\}$, $\psi_j$ is a weight of $N_{\phi_j}$.

It follows that both $P_{\phi}/\iota^{-1}(X)$ and $S_{\phi_1}\otimes \dots \otimes S_{\phi_n}$ have 
the same set of weights. Since, in addition, these modules have one
dimensional weight spaces, the surjective map must be an isomorphism. But, as quoted above, 
$S_{\phi_1}\otimes \dots \otimes S_{\phi_n}$ is simple, which implies 
$N_{\phi} = \iota^{-1}(X)$. This proves the statement.\qed
\end{enumerate}

\section{Quantum Weyl algebras} \label{section-qwa}

In this section we will consider two versions of quantum Weyl algebras. 
For a short account on the description of these quantum Weyl algebras as ring of skew differential operators, see Section \ref{classification}.
In the first subsection, we recall Maltsiniotis' quantum Weyl algebra. In the second we introduce 
the version of Akhavizadegan and Jordan. In the last subsection, we recall the isomorphism between 
the convenient localisations of these algebras.  

For this, we need to introduce quantum integers. The reader is referred to \cite{AJ} and \cite{J}
for a detailed exposition. Some complements may be found in \cite{R1}. \\

Let $q\in\k^\ast$. 
For any integer $i\in\N^\ast$, we let
\[
(i)_q = 1 +q + \dots + q^{i-1}, 
\]
and extend this definition by $(0)_q=0$. Further, for $j\in\Z\setminus\N$, we set
\[
(j)_q = -q^j(-j)_q = -\left(q^j + q^{j+1} + \dots q^{-1}\right). 
\]
These scalars will be called {\em quantum integers associated to} $q$ or, for simplicity, $q$-{\em integers}.

\begin{remark} -- \rm
\begin{enumerate}
\item If $q=1$, then, for all $i\in\Z$, $(i)_q=i$.  
\item If $q\neq 1$, then,
\[
\forall i\in\Z \quad\quad (i)_q = \displaystyle\frac{1-q^i}{1-q} . 
\]
\end{enumerate}
\end{remark}

Recall that an $n \times n$ matrix $(\lambda_{ij})$ with entries in $\k^\ast$ is called 
skew-symmetric provided $\lambda_{ii}=1$ and $\lambda_{ij}=\lambda_{ji}^{-1}$, 
for all $1 \le i,j \le n$.

\subsection{The quantum Weyl algebra of G.~Maltsiniotis}

\begin{subdefinition} -- Let $n\in\N^\ast$. 
To any $n$-tuple $\bar{q}=(q_{1} ,...,q_n)\in(\k^\ast)^n$ and any $n \times n$ skew symmetric matrix 
$\Lambda=(\lambda_{ij}) \in M_{n}(k)$, we associate the $\k$-algebra
$A_n^{\bar{q},\Lambda}$ with generators $y_1,...,y_n,x_1 ,...,x_n$ and relations:
\begin{equation}\label{awqm-1}
\begin{array}{lll}
\forall (i,j)\in\N^2, \; 1 \le i<j\le n,  & 
x_{i} x_{j}=\lambda _{ij} q_{i} x_{j} x_{i},  &  y_{i}  y_{j}=\lambda _{ij} y_{j} y_{i}, \cr   
&  x_{i} y_{j}=\lambda _{ij}^{-1} y_{j} x_{i},  & y_{i} x_{j}=\lambda _{ij}^{-1} q_{i}^{-1} x_{j} y_{i}. 
\end{array}
\end{equation}
\begin{equation}\label{awqm-2}
\begin{array}{l}
\forall i\in\N, \; 1 \le i \le n,\;  x_{i} y_{i} -q_{i}y_{i} x_{i} =1+\sum 
_{j=1}^{i-1} (q_{j}-1)y_{j} x_{j}.
\end{array}
\end{equation}
Further, we put $z_0=1$ and, for $1 \le i \le n$,  
\[
z_i =x_iy_i - y_ix_i. 
\]
\end{subdefinition}

The family $\{z_1,\dots,z_n\}$ will play a central role in computations.

\begin{subremark} -- \rm It is easy to see that for all $i$, $1\le i \le n$:
\begin{equation}
\begin{array}{l}
\qquad z_{i}=1+\sum_{j=1}^{i} (q_{j}-1)y_{j} x_{j} 
\qquad\mbox{and}\qquad z_{i}=z_{i-1}+(q_{i}-1)y_{i}x_{i}.
\end{array}
\end{equation}
As a consequence, relations (\ref{awqm-2}) read: \\ 
\begin{equation}
\begin{array}{l}
\forall i, 1\le i \le n, \;\; x_{i} y_{i} -q_{i} y_{i} x_{i} =z_{i-1}
\end{array}
\end{equation}
\end{subremark}

\begin{subproposition} -- \label{maltsiniotis}
The algebra $A_n^{\bar{q},\Lambda}$ is a noetherian integral domain. In addition, the family
$\{y_{1}^{j_{1}}x_{1}^{i_{1}}...y_{n}^{j_{n}}x_{n}^{i_{n}},\; (j_{1},i_{1},...,j_{n},i_{n}) \in \N^{2n}\}$
is a 
PBW basis of $A_n^{\bar{q},\Lambda}$.
\end{subproposition}

\proof See \cite[Remarque 2.1.3]{R1} and \cite[1.9]{AJ}. \qed

\begin{subproposition}\label{maltsiniotis1} -- 
For all $1 \le i \le n$, $z_{i}$ is a normal element of $A_n^{\bar{q},\Lambda}$. More precisely, we have:
\begin{itemize}
\item 
for all $i, j$, $1 \le i  <  j \le n, \qquad z_{i}x_{j}=x_{j}z_{i} \qquad \mbox{and}\qquad  z_{i}y_{j}=y_{j}z_{i},$
\item 
for all $i, j$, $1 \le j \le i \le n, \qquad z_{i}x_{j}=q_{j}^{-1}x_{j}z_{i} \qquad \mbox{and} \qquad z_{i}y_{j}=q_{j}y_{j}z_{i},$
\item 
for all $i, j$, $1 \le i,j \le n, \qquad z_{i}z_{j}=z_{j}z_{i}.$
\end{itemize}
\end{subproposition}

\proof See \cite[p.~285]{J}. \qed\\

As pointed out above, the multiplicative set generated by the elements $z_i$, $1 \le i \le n$, consists of normal elements.
It follows that it is an Ore set and that we may form the corresponding localisation.

\begin{subdefinition} -- 
We denote by $B_{n}^{\bar{q},\Lambda}$ the localisation of $A_{n}^{\bar{q}, \Lambda}$ at the multiplicative set generated
by the elements $z_i$, $1 \le i \le n$.
\end{subdefinition}

Since $A_{n}^{\bar{q},\Lambda}$ is an integral domain, there is a canonical injective morphism of $\k$-algebras 
\[
A_{n}^{\bar{q},\Lambda} \stackrel{can.inj.}{\longrightarrow} B_{n}^{\bar{q},\Lambda}.
\]
In the sequel, we will often identify an element of 
$A_{n}^{\bar{q},\Lambda}$ with its image in $B_{n}^{\bar{q},\Lambda}$ under the above canonical injection.

Denote by $R$ (resp $R^\circ$) the $\k$-subalgebra of $A_n^{\bar{q},\Lambda}$ 
(resp. $B_n^{\bar{q},\Lambda}$) generated by $z_1,\dots,z_n$ (resp. $z_1,\dots,z_n$ and their inverses). Hence, $R$ and $R^\circ$ are commutative $\k$-algebras. 

\begin{subproposition}\label{maltsiniotis2} -- Assume $q_i\neq 1$ for all $1 \le i \le n$.
The elements $z_1,\dots,z_n$ are algebraically independent in $R$. Moreover, the set  
\[
\b=\left\{y_{1}^{j_{1}}x_{1}^{i_{1}}...y_{n}^{j_{n}}x_{n}^{i_{n}},\; (j_{1},i_{1},...,j_{n},i_{n}) \in \N^{2n}, \;
i_kj_k=0, \, \forall 1 \le k \le n\right\} 
\]
is a basis of $A_n^{\bar{q},\Lambda}$ regarded as a left or as a right $R$-module.
Further, $\b$ is a basis of $B_n^{\bar{q},\Lambda}$ both as a left and as a right $R^\circ$-module.
\end{subproposition}

\proof
We leave the detailed proof to the interested reader. All the statements follow from a careful use of the above relations and Proposition \ref{maltsiniotis}. \qed

\subsection{The quantum Weyl algebra of M. Akhavizadegan \& D. Jordan}

\begin{subdefinition} -- Let $n\in\N^\ast$. 
To any $n$-tuple $\bar{q}=(q_{1} ,...,q_n)\in(\k^\ast)^n$ and any $n \times n$ skew symmetric matrix 
$\Lambda=(\lambda_{ij}) \in M_{n}(k)$, we associate the $\k$-algebra
$\A_n^{\bar{q},\Lambda}$ with generators $y_1,...,y_n,x_1 ,...,x_n$ and relations:
\begin{equation}\label{awqmAJ-1}
\begin{array}{lll}
\forall (i,j)\in\N^2, \; 1 \le i<j\le n,  & 
x_{i} x_{j}=\lambda _{ij} x_{j} x_{i},  &  y_{i}  y_{j}=\lambda _{ij} y_{j} y_{i}, \cr   
&  x_{i} y_{j}=\lambda _{ij}^{-1} y_{j} x_{i},  & y_{i} x_{j}=\lambda _{ij}^{-1} x_{j} y_{i}. 
\end{array}
\end{equation}
\begin{equation}\label{awqmAJ-2}
\begin{array}{l}
\forall i\in\N, \; 1 \le i \le n,\;  x_{i} y_{i} -q_{i}y_{i} x_{i} = 1.
\end{array}
\end{equation}
Let us also define $z_0=1$ and, for $i$ such that $1 \le i \le n$,  
\[
z_i =x_iy_i - y_ix_i. 
\]
\end{subdefinition}

As before, these elements will play a crucial role.

\begin{subremark} -- \rm \label{rm4relations}
It is easy to see that the four relations in (\ref{awqmAJ-1}) actually hold 
for all $(i,j)\in\N^2$ such that $1 \le i \neq j \le n$ and that:
\[
\forall i, \qquad 1 \le i \le n, \qquad z_{i}=1+(q_i-1)y_i x_i .
\]
\end{subremark}

Two results analogous to Propositions \ref{maltsiniotis1} and \ref{maltsiniotis2} are also true in this case.

\begin{subproposition} -- The algebra $\A_n^{\bar{q},\Lambda}$ is a noetherian integral domain. In addition, the family
$\{y_{1}^{j_{1}}x_{1}^{i_{1}}...y_{n}^{j_{n}}x_{n}^{i_{n}},\; (j_{1},i_{1},...,j_{n},i_{n}) \in \N^{2n}\}$
is a $\k$-basis of $\A_n^{\bar{q},\Lambda}$.
\end{subproposition}

\proof See \cite[1.9]{AJ}. \qed\\

In particular, the family $\{z_1,\dots, 
z_n\}$ verifies properties similar to those of the previous subsection.

\begin{subproposition}\label{z-et-xy} -- 
For all $1 \le i \le n$, $z_{i}$ is a normal element of $\A_n^{\bar{q},\Lambda}$. More precisely, we have:
\begin{itemize}
\item 
for all $i, j$, $1 \le i \neq  j \le n, \qquad z_{i}x_{j}=x_{j}z_{i} \qquad \mbox{and}\qquad  z_{i}y_{j}=y_{j}z_{i},$
\item 
for all $i$, $1 \le i \le n, \qquad z_{i}x_{i}=q_{i}^{-1}x_{i}z_{i} \qquad \mbox{and} \qquad z_{i}y_{i}=q_{i}y_{i}z_{i},$
\item 
for all $i, j$, $1 \le i,j \le n, \qquad z_{i}z_{j}=z_{j}z_{i}.$
\end{itemize}
\end{subproposition}

\proof See \cite[p. 286]{AJ}. \qed\\

As pointed out before, the multiplicative set generated by all the $z_i$'s consists of normal elements.
Hence, it is an Ore set and thus we may form the corresponding localisation.

\begin{subdefinition} -- 
We denote by $\B_{n}^{\bar{q},\Lambda}$ the localisation of $\A_{n}^{\bar{q}, \Lambda}$ at the multiplicative set generated
by the elements $z_i$, $1 \le i \le n$.
\end{subdefinition}

Notice that, again, the fact that $\A_{n}^{\bar{q},\Lambda}$ is an integral domain implies that
the canonical morphism of $\k$-algebras 
\[
\A_{n}^{\bar{q},\Lambda} \stackrel{can.inj.}{\longrightarrow} \B_{n}^{\bar{q},\Lambda} 
\]
is injective, and we will identify any element of 
$\A_{n}^{\bar{q},\Lambda}$ with its image in $\B_{n}^{\bar{q},\Lambda}$ under this injection.

Denote by $R$ (resp $R^\circ$) the $\k$-subalgebra of $\A_n^{\bar{q},\Lambda}$ 
(resp. $\B_n^{\bar{q},\Lambda}$) generated by $z_1,\dots,z_n$ (resp. and their inverses). 

\begin{subproposition} -- \label{R-basis-of-A}
Assume $q_i\neq 1$ for all $1 \le i \le n$.
The elements $z_1,\dots,z_n$ are algebraically independent in $R$. Further, the set
\[
\b=\{y_{1}^{j_{1}}x_{1}^{i_{1}}...y_{n}^{j_{n}}x_{n}^{i_{n}},\; (j_{1},i_{1},...,j_{n},i_{n}) \in \N^{2n}, \;
i_kj_k=0, \, \forall 1 \le k \le n\}. 
\]
is a basis of $\A_n^{\bar{q},\Lambda}$ regarded as a left and as a right $R$-module.
Moreover, $\b$ is also a basis of $\B_n^{\bar{q},\Lambda}$ as a left and as a right $R^\circ$-module.
\end{subproposition}

\proof 
It is similar to the proof of Prop. \ref{maltsiniotis2}.
\qed

\subsection{A common localisation for quantum Weyl algebras}

We finish this section by recalling a result which states that the localisations of both versions of quantum Weyl algebras are isomorphic.

\begin{subtheorem} -- Let $n\in\N^\ast$. 
For any $n$-tuple $\bar{q}=(q_{1} ,...,q_n)\in(\k^\ast)^n$ and any $n \times n$ skew symmetric matrix 
$\Lambda=(\lambda_{ij}) \in M_{n}(\k)$, the map defined by:
\[
\begin{array}{ccrcl}
\theta & : & \B_n^{\bar{q},\Lambda} & \longrightarrow & B_n^{\bar{q},\Lambda} \cr
 & & y_i & \mapsto & y_i \cr
 & & x_i & \mapsto & z_{i-1}^{-1}x_i \cr
 & & z_i & \mapsto & z_{i-1}^{-1}z_i
\end{array} 
\]
is an isomorphism of $\k$-algebras.

\end{subtheorem}

\proof See p. 287 of \cite{AJ}.\qed

\section{Weight modules for quantum Weyl algebras}

In this section, we fix an $n$-tuple $\bar{q}=(q_{1} ,...,q_n)\in(\k^\ast)^n$ and an 
$n \times n$ skew symmetric matrix 
$\Lambda=(\lambda_{ij}) \in M_{n}(\k)$ and we associate to this data the algebra 
$\B_n^{\bar{q},\Lambda}$ introduced in Section \ref{section-qwa}.
We consider the context of Section \ref{section-weight-categories} for the pair 
$(R^\circ,\B_n^{\bar{q},\Lambda})$, where $R^\circ$ is the $\k$-subalgebra of 
$\B_n^{\bar{q},\Lambda}$ generated by the $z_i$'s and their inverses. 
As a consequence of Proposition \ref{R-basis-of-A}, $R^\circ$
is a commutative Laurent polynomial $\k$-algebra in the indeterminates $z_1,\dots,z_n$. 

As we will see, the situation will heavily depend on whether some $q_i$ is a root of unity in $\k$
or not. We will say that the $n$-tuple $\bar{q}$ is {\em generic} whenever, 
for all integer $i$, $1 \le i \le n$, $q_i$ is not a root of unity. 

In the first subsection, we show that the pair $(R^\circ,\B_n^{\bar{q},\Lambda})$ is indeed a 
strongly-free extension. We also consider the action of the canonical generators of 
$\B_n^{\bar{q},\Lambda}$ on the basis attached to this extension. From that basis arises in a 
natural way an action of the group $\Z^n$ on the weights. We investigate this action in the
second subsection and show for example that it is free if and only if no $q_i$ is a root of unity.
As we discussed in the first section, each character $\phi$ of $R^\circ$ gives rise to a projective 
weight $\B_n^{\bar{q},\Lambda}$-module $P_\phi$. The structure of these modules is the subject of 
the last subsection. We investigate their weights in connection with the above group action, the 
dimension of their weight spaces, characterise their simplicity and start to discuss whether two such 
modules are isomorphic or not. 

\subsection{Basic results on the pair $(R^\circ,\B_n^{\bar{q},\Lambda})$}

By Proposition \ref{R-basis-of-A}, 
\[
\b=\left\{y_{1}^{j_{1}}x_{1}^{i_{1}}...y_{n}^{j_{n}}x_{n}^{i_{n}},\; (j_{1},i_{1},...,j_{n},i_{n}) 
\in \N^{2n}, \;
i_kj_k=0, \, \forall k, 1 \le k \le n\right\}. 
\]
is a basis of $\B_n^{\bar{q},\Lambda}$ regarded either as a left or as a right $R^\circ$-module. It 
will be convenient 
to write this basis in a slightly different way. Let $k=(k_1,\dots,k_n)\in\Z^n$. We denote
\[
b_k = \prod_{1 \le i \le n}  b_{k,i}
\qquad
\mbox{where, for $1 \le i \le n$}, 
\qquad
b_{k,i} 
=
\left\{
\begin{array}{ccc}
x_i^{k_i} & \mbox{if} & k_i \ge 0, \cr
y_i^{-k_i} & \mbox{if} & k_i \le 0 . 
\end{array}
\right.
\]
With this notation, we have
\[
\b=\{b_k, \, k\in\Z^n\}. 
\]
Given $k\in\Z^n$, let $\sigma_k$ be the $\k$-algebra automorphism of $R^\circ$ such that
\[
\forall i, \ 1 \le i \le n, \quad\sigma_k(z_i) = q_i^{-k_i}z_i. 
\]
From relations in Proposition \ref{z-et-xy}, we get that
\[
\forall\, k\in\Z^n, \quad \forall r\in R^\circ, \quad rb_k=b_k\sigma_k(r).  
\]
In the notation of Definition \ref{definition-strongly-free}, this means that, for all 
$k\in\Z^n$, $\sigma_k=\sigma_{b_k}$. 
These facts afford the following lemma.

\begin{sublemma} -- 
The pair $(R^\circ,\B_n^{\bar{q},\Lambda})$ is of strongly-free type.
\end{sublemma}

We finish this subsection by examining how multiplication by generators modify the elements of $\b$.
From now on, we let $\{e_1,\dots,e_n\}$ stand for the canonical basis of $\Z^n$. 

\begin{sublemma} -- \label{relations-x-y-et-b} 
Assume $q_i\neq 1$ for all $1 \le i \le n$.
Given $k\in\Z^n$ and $i$ such that $1 \le i \le n$,
\[
x_ib_k = 
\left\{
\begin{array}{lll}
\left(\displaystyle\prod_{1 \le j < i} \lambda_{ij}^{k_j} \right) b_{k+e_i} & \mbox{if} & k_i \ge 0, \cr
\left(\displaystyle\prod_{1 \le j < i} \lambda_{ij}^{k_j} \right) b_{k+e_i}\displaystyle\frac{q_i^{-k_i}z_i-1}{q_i-1} & \mbox{if} & k_i < 0 \cr
\end{array}
\right. 
\]
and
\[
y_ib_k = 
\left\{
\begin{array}{lll}
\left(\displaystyle\prod_{1 \le j < i} \lambda_{ij}^{k_j} \right)^{-1} b_{k-e_i}\displaystyle\frac{q_i^{-k_i+1}z_i-1}{q_i-1} & \mbox{if} & k_i > 0, \cr 
\left(\displaystyle\prod_{1 \le j < i} \lambda_{ij}^{k_j} \right)^{-1} b_{k-e_i} & \mbox{if} & k_i \le 0.
\end{array}
\right.
\]
\end{sublemma}

\proof This is straightforward using relation (\ref{awqmAJ-1}), (\ref{awqmAJ-2}) and the equalities 
of Remark \ref{rm4relations} and Proposition \ref{z-et-xy}.\qed

\subsection{Action of $\Z^n$ on weights arising from the basis $\b$}\label{action-B}

Recall from Remark \ref{action-automorphismes-sur-poids} the right action of $\Aut_\k(R^\circ)$ on $\widehat{R^\circ}$.
On the other hand, we have group morphisms as follows
\[
\epsilon: \Z^n \rightarrow  (\k^\ast)^n \qquad \quad (k_1,\dots,k_n) \mapsto (q_1^{-k_1},\dots,q_n^{-k_n})
\]
and
\[
\delta: (\k^\ast)^n \rightarrow \Aut_\k(R^\circ) \qquad (\lambda_1,\dots,\lambda_n) \mapsto\left( z_i \mapsto \lambda_i z_i \right)_{1 \le i \le n}.
\]
\smallskip

Now, composing these maps gives rise to 
a left action of $(\k^\ast)^n$ on $\widehat{R^\circ}$ and 
a left action of $\Z^n$ on $\widehat{R^\circ}$ given by
\[
(\k^\ast)^n 
\stackrel{\delta}{\longrightarrow} \Aut_\k(R^\circ)
\stackrel{}{\longrightarrow} \S(\widehat{R^\circ})
\qquad\mbox{and}\qquad
\Z^n 
\stackrel{\epsilon}{\longrightarrow} (\k^\ast)^n 
\stackrel{\delta}{\longrightarrow} \Aut_\k(R^\circ)
\stackrel{}{\longrightarrow} \S(\widehat{R^\circ}),
\]
where $\S(\widehat{R^\circ})$ is the group of permutations of the set $\widehat{R^\circ}$. Alternatively, the latter
may be described as follows:
\begin{equation}
\begin{array}{ccrcl}
 & & \Z^n \times \widehat{R^\circ} & \longrightarrow & \widehat{R^\circ} \cr 
 & & (k,\phi) & \mapsto & \phi \circ \sigma_k.
\end{array}
\end{equation}

\smallskip
These actions verify some easily proven properties that we list in the following lemma.

\begin{sublemma} --
Using the notation above, the following holds:
\begin{enumerate}
\item
the left action of $(\k^\ast)^n$ on $\widehat{R^\circ}$ is faithful; 
\item the morphism $\epsilon  \, : \Z^n \longrightarrow (\k^\ast)^n$ is injective if and only if  $\bar{q}$ is generic;
\item the left action of $\Z^n$ on $\widehat{R^\circ}$ is faithful if and only if 
$\bar{q}$ is generic.
\end{enumerate}
\end{sublemma}

Fixing $\phi\in \widehat{R^\circ}$, we obtain a map $\Z^n \to \widehat{R^\circ}$.

\begin{sublemma} -- \label{genericity-injectivity}
Let $\phi\in\widehat{R^\circ}$. The map $\Z^n \to \Z^n.\phi\subseteq\widehat{R^\circ}$ 
sending $k$ to $\phi\circ\sigma_k$ is injective if and only if $\bar{q}$ is generic. 
\end{sublemma}

\proof Let $k\in\Z^n$. By definition, given $i$ such that $1 \le i \le n$, we know that 
$\phi\circ\sigma_k(z_i)=q_i^{-k_i}\phi(z_i)$. 
But, $z_i$ being an invertible element of $\B_n^{\bar{q},\Lambda}$, $\phi(z_i)$ is necessarily non zero. The result follows. \qed

\begin{subcorollary} -- \label{generic-and-free}
The following conditions are equivalent:
\begin{enumerate}
\item
$\bar{q}$ is generic; 
\item there exists $\phi\in\widehat{R^\circ}$ such that $\stab_{\Z^n}(\phi)=\{0_{\Z^n}\}$;
\item for all $\phi\in\widehat{R^\circ}$, $\stab_{\Z^n}(\phi)=\{0_{\Z^n}\}$, in other words
the action of $\Z^n$ on $\widehat{R^\circ}$ is free.
\end{enumerate}
\end{subcorollary}

\proof This is an immediate consequence of Lemma \ref{genericity-injectivity}.\qed

\subsection{Structure of $P_\phi$ for $\phi\in\widehat{R^\circ}$}\label{base-de-P}

Let $\phi\in\widehat{R^\circ}$. 
Recall that $P_\phi = \B_n^{\bar{q},\Lambda} \otimes_{R^\circ} \k_\phi$, where 
$\k_\phi=R^\circ/\ker(\phi)$. Given $k\in\Z^n$, set
\[
v_k = b_k \otimes_R \ov{1} .    
\]
By Remark \ref{structure-P}, the set $\{v_k,\, k\in\Z^n\}$, 
is a $\k$-basis of $P_\phi$ and 
\begin{equation}\label{inclusion-weight}
\forall k\in\Z^n, \quad \k v_k \subseteq (P_\phi)(\phi\circ\sigma_k) . 
\end{equation}
Hence, $P_\phi$ is a weight module of $\B_n^{\bar{q},\Lambda}$ and the set of weights of $P_\phi$
is the $\Z^n$-orbit of $\phi$, that is $\Z^n.\phi$. 

\begin{subproposition} -- \label{generic-weight-decomposition}
Let $\phi\in\widehat{R^\circ}$. If $\bar{q}$ is generic, then,
for all $k\in\Z^n, \quad \k v_k = (P_\phi)(\phi\circ\sigma_k)$. In particular, 
the weight decomposition of the weight $\B_n^{\bar{q},\Lambda}$-module $P_\phi$ is:
\[
P_\phi = \bigoplus_{k\in\Z^n} \k v_k = \bigoplus_{k\in\Z^n} (P_\phi)(\phi\circ\sigma_k).  
\]
\end{subproposition}

\proof Given $k\in\Z^n$, we know from the previous discussion that 
$\k v_k \subseteq (P_\phi)(\phi\circ\sigma_k)$. 
The result is thus an obvious consequence of Lemma \ref{genericity-injectivity}. \qed

\begin{subremark} -- \rm \label{generic-P-weyl} We may actually strengthened 
Proposition \ref{generic-weight-decomposition} as follows. 
Let $\phi\in\widehat{R^\circ}$. The following statements are equivalent:\\
(i) $\bar{q}$ is generic; \\
(ii) the map $\b \mapsto \widehat{R^\circ}$, $b \mapsto \phi \circ \sigma_b$ is injective;\\
(ii) the weight spaces of $P_\phi$ all have dimension at most one.\\
This follows from the above discussion taking Lemma \ref{genericity-injectivity} and the inclusion 
(\ref{inclusion-weight}) into account.
\end{subremark}

\begin{subremark} -- \rm \label{N-phi}
Supposing that $\bar{q}$ is generic, the hypotheses of Remark \ref{generic-P} are fulfilled.
Hence, we know that, for all $\phi\in\widehat{R^\circ}$, the module $P_\phi$
has a maximum strict weight $\B_n^{\bar{q},\Lambda}$-submodule $N_\phi$ and a unique simple quotient 
$S_\phi=P_\phi/N_\phi$. 
\end{subremark}

We continue by giving explicit expressions for the action of $x_i$ and $y_i$ on basis vectors of 
$P_\phi$. 
Recall that $\{e_1,\dots,e_n\}$ is the canonical basis of $\Z^n$.

\begin{sublemma} -- \label{action-xy}
Assume $q_i\neq 1$ for all $1 \le i \le n$.
Let $\phi\in\widehat{R}$. Let $1 \le i \le n$ and let $k\in\Z^n$~: 
\begin{enumerate}
\item
the action of $x_i$ on $v_k$ is given by~:
\begin{equation} \label{action-x}
x_i . v_k 
=
\left\{
\begin{array}{rcc}
\left(\displaystyle\prod_{1 \le j < i} \lambda_{ij}^{k_j} \right)v_{k+e_i} & \mbox{if} & k_i \ge 0, \cr
\displaystyle\frac{q_i^{-k_i}\phi(z_i)-1}{q_i-1} \left(\displaystyle\prod_{1 \le j < i} \lambda_{ij}^{k_j} \right)v_{k+e_i} & \mbox{if} & k_i < 0 ; 
\end{array}
\right.
\end{equation} 
\item the action of $y_i$ on $v_k$ is given by~:
\begin{equation} \label{action-y}
y_i . v_k 
= \left\{
\begin{array}{rcc}
\displaystyle\frac{q_i^{-k_i+1}\phi(z_i)-1}{q_i-1} \left(\displaystyle\prod_{1 \le j < i} \lambda_{ij}^{k_j} \right)^{-1}v_{k-e_i} & \mbox{if} & k_i > 0, \cr
\left(\displaystyle\prod_{1 \le j < i} \lambda_{ij}^{k_j} \right)^{-1}v_{k-e_i} & \mbox{if} & k_i \le 0. 
\end{array} 
\right.
\end{equation}
\end{enumerate}
\end{sublemma}

\proof It is an immediate consequence of Lemma \ref{relations-x-y-et-b}. \qed\\

It is clear after Lemma \ref{action-xy} that a key feature of $P_\phi$ depends on whether the scalars appearing in relations 
(\ref{action-x}) and (\ref{action-y}) may vanish or not, which motivates the next definition.

\begin{subdefinition} -- Let $\phi\in\widehat{R^\circ}$. We define the complexity of $\phi$ as the set
\[
\comp(\phi) = \{i \in\{1,\dots,n\} \tq \phi(z_i) \in \la q_i \ra\}.
\]
\end{subdefinition}

\begin{subremark} -- \rm
Any two elements of $\widehat{R^\circ}$ in the same $\Z^n$-orbit have the same complexity. 
\end{subremark}

Recall that, after Proposition \ref{sous-objets-et-poids}, any submodule of $P_\phi$ is a weight submodule.
The next proposition shows that weights with empty complexity give rise to simple modules.

\begin{subproposition} -- \label{objets-de-niveau-0} 
Suppose that $\bar{q}$ is generic.
If $\phi\in\widehat{R^\circ}$, the following are equivalent: 
\begin{enumerate}
\item
the complexity of $\phi$ is empty, 
\item $P_\phi$ is a simple object both of $\B_n^{\bar{q},\Lambda}-\Mod$ and of $\B_n^{\bar{q},\Lambda}-\wMod$. 
\end{enumerate}
\end{subproposition}
\proof By Proposition \ref{sous-objets-et-poids}, it is enough to consider the simplicity of $P_\phi$
as an object of $\B_n^{\bar{q},\Lambda}-\wMod$. 

Suppose $M$ is a non zero subobject of $P_\phi$ in the category 
$\B_n^{\bar{q},\Lambda}-\wMod$. By Proposition 
\ref{generic-weight-decomposition}, it must contain one of the $v_k$'s with $k\in\Z^n$. 
This in turn implies, using Lemma \ref{action-xy}, that if the complexity of $\phi$ is empty, then $v_0\in M$ and thus $M=P_\phi$. 

Conversely, suppose the complexity of $\phi$ is not empty and let $i$ be an integer belonging to $\comp(\phi)$. 
Thus, $\phi(z_i) = q_i^\ell$ for some $\ell\in\Z$. 
Suppose first that $\ell < 0$. It follows from Lemma \ref{action-xy} that $x_i.v_k=0$ whenever $k_i=\ell$.
As a consequence, $\oplus_{k, k_i \le \ell} \k v_k$ is a strict $\B_n^{\bar{q},\Lambda}$-submodule of $P_\phi$.
Suppose now that $\ell \ge 0$. In this case, it follows from Lemma \ref{action-xy} that 
$y_i.v_k=0$ whenever $k_i=\ell+1$.
Consequently, $\oplus_{k, k_i \ge \ell+1} \k v_k$ is a strict $\B_n^{\bar{q},\Lambda}$-submodule of $P_\phi$.
Hence, in all cases, $P_\phi$ is not simple. \qed\\

In the next statement, we analyse the link between $P_\phi$ and $P_{\phi\circ\sigma_{e_\ell}}$ for an integer $\ell$ such that
$1 \le \ell \le n$. Recall that $e_i$ stands for the $i$-th vector of the canonical basis of $\Z^n$.  

\begin{subproposition} -- \label{iso-des-P-phi} 
Suppose $q_i\neq 1$ for all $1 \le i \le n$. 
Fix an integer $\ell$ such that $1 \le \ell \le n$
and $\phi\in\widehat{R^\circ}$. 
\begin{enumerate}
\item
If $(\lambda_k)_{k\in\Z^n} \in \k^{(\Z^n)}$ satisfies the following conditions: 
\begin{enumerate}
\item
$\lambda_k = \lambda_{k+e_\ell}$, if $k_\ell \neq -1$;
\item $\lambda_k = \displaystyle\frac{\phi(z_\ell)-1}{q_\ell-1}\lambda_{k+e_\ell}$, if $k_\ell = -1$;
\item $\lambda_k = \lambda_{k+e_i}$, if $i<\ell$;
\item $\lambda_k = \lambda_{i\ell}^{-1}\lambda_{k+e_i}$, if $i>\ell$,
\end{enumerate}
then, the $\k$-linear map
\[
\begin{array}{ccrcl}
 & & P_{e_\ell.\phi} & \longrightarrow & P_\phi \cr
 & & w_k & \mapsto & \lambda_kv_{k+e_\ell},
\end{array}
\]
where $(w_k)_{k\in\Z^n}$ and $(v_k)_{k\in\Z^n}$ are respectively the canonical bases of 
$P_{e_\ell.\phi}$ and $P_\phi$ --as defined at the beginning of this subsection-- is a morphism in 
$\B_n^{\bar{q},\Lambda}-\wMod$.
\item Assume in addition that $\phi(z_\ell) \neq 1$. Then,  
$P_{e_\ell.\phi}$ is isomorphic to $P_\phi$ in $\B_n^{\bar{q},\Lambda}-\wMod$.
\end{enumerate}
\end{subproposition}

\proof The first item is a lengthy straightforward verification.
The second item follows directly from the first one.
\qed

\section{Examples: the case where all the entries of $\Lambda$ equal $1$} \label{section-exemples}

As it will be shown in Section \ref{QWa-et-ZT}, it turns out that the study of simple weight modules 
of the extension $(R^\circ,\B_n^{\bar{q},\Lambda})$ may be reduced to the case where the skew-
symmetric matrix $\Lambda$ has all its entries equal to $1$. For this reason we concentrate, in the
present section, to this special case. It is easier to handle since under this hypothesis, 
$\B_n^{\bar{q},\Lambda}$ is a tensor product of $n$ copies of $\B_1^{\bar{q},\Lambda}$. 

Hence, we first investigate the case where $n=1$. We give a complete and explicit description of
the projective weight modules $P_\phi$ attached to each character $\phi$ of $R^\circ$. From this we 
deduce a complete description of the simple weight modules and classify them. The second subsection 
comes back to the case of an arbitrary integer $n$. It is shown that simple weight modules in this 
case may be explicitly described using simple weight modules arising in the case $n=1$. This is illustrated in the third subsection where the case $n=2$ is studied in full details using very intuitive graphs to describe the modules $P_\phi$.

In this section, we suppose that $\bar{q}$ is generic. We also assume that all the entries of the 
skew-symmetric matrix
$\Lambda$ equal $1$.
Under this last hypothesis, there is an isomorphism of $\k$-algebras
\[
\begin{array}{ccrcl}
& & \B_n^{\bar{q},\Lambda} & \longrightarrow & \B_1^{q_1} \otimes\dots\otimes \B_1^{q_n}  
\end{array}
\]
that sends each $x_i$ to the elementary tensor 
$1 \otimes\dots\otimes x_1 \otimes\dots\otimes 1$, 
and each $y_i$ to $1 \otimes\dots\otimes y_1 \otimes\dots\otimes 1$, 
where $x_1$ and $y_1$ are, respectively, in the $i$-th position, where, for $q\in\k^\ast$, we put 
$\B_1^q=\B_1^{\bar{q},\Lambda}$ with $\bar{q}=(q)$ and $\Lambda=(1)$ .
This isomorphism clearly induces in turn the isomorphism
\[
\begin{array}{ccrcl}
& & \k[z_1^{\pm 1},\dots,z_n^{\pm 1}] & \longrightarrow & \k[z_1^{\pm 1}] \otimes\dots\otimes \k[z_1^{\pm 1}]. 
\end{array}
\]
We are now ready to apply the results of Proposition \ref{weight-in-tpa} and Subsection \ref{TP-of-SFE}.

\medskip

From now on, 
we denote by $\1$ the unique element of $\widehat{\k[z_1^{\pm 1},\dots,z_n^{\pm 1}]}$ that sends $z_1,\dots,z_n$ 
to $1$. 

\subsection{The case $n=1$}\label{exemple-n=1}

The observations in the Introduction show that one can build simple weight modules
of $\B_n^{\bar{q},\Lambda}$ from those appearing  in the case $n=1$. So we start by studying this case.
In this subsection, we suppose $n=1$.

\begin{subremark} -- \rm \label{complexite-orbite-1}
Let $\phi\in \widehat{R^\circ}$. 
It is clear that the following statements are equivalent: 
\begin{enumerate}
\item
$\comp(\phi) = \{1\}$; 
\item
 $\1 \in \Z.\phi$;
\item
 $\Z.\phi=\Z.\1$.
\end{enumerate}
\end{subremark}

Fix now an element $\phi\in\widehat{R^\circ}$. 
As a consequence of Proposition \ref{generic-weight-decomposition} and 
Corollary \ref{generic-and-free}, 
the set of weights of $P_\phi$ is $\Z.\phi$, with $\Z.\phi$ being equipotent to $\Z$  
and each weight space of $P_\phi$ is one dimensional over $\k$. \\

According to Proposition \ref{objets-de-niveau-0}, $P_\phi$ is simple if and only if 
the complexity of $\phi$ is empty, that is if and only if $\phi$ is in the orbit of $\1$ under the
action of $\Z$ (cf. Remark \ref{complexite-orbite-1}). We thus analyse these two cases. 

In order to 
visualize the action under consideration, we represent it by a graph as follows.
The vectors of the canonical basis of $P_\phi$ are displayed horizontally along a line
and constitute the set of vertices of the graph. 
By Lemma \ref{action-xy}, $x_1$ (resp. $y_1$) acts on a given $v_\alpha$ sending it to 
a scalar multiple of $v_{\alpha +1}$ (resp. $v_{\alpha -1}$). In case the corresponding scalar is 
nonzero, we draw an arrow, labeled by the corresponding generator, from $v_\alpha$ to
$v_{\alpha +1}$ (resp. $v_{\alpha -1}$); in case it is zero, we do not draw any arrow.  
A careful observation allows to "see" in this graph the unique maximal submodule and corresponding 
simple quotient. 

\begin{subexample} -- \label{exemple-hors-1} \rm
Suppose $\phi\not\in\Z.\1$ or, equivalently, $\comp(\phi)=\emptyset$. 
\begin{enumerate}
\item By Proposition \ref{objets-de-niveau-0}, $P_\phi$ is a simple $\B_1^{\bar{q},\Lambda}$-weight module. Its corresponding scheme is as follows:
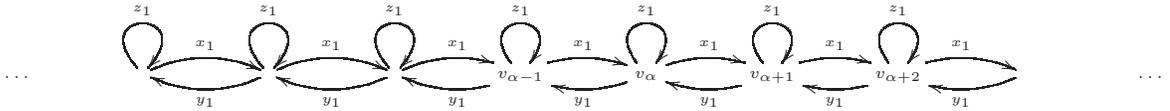
\tiny
\[
\xymatrix@!C{
\dots
& \ar@(ul,ur)^{z_1} \ar@/^/[r]^{x_1}
& \ar@(ul,ur)^{z_1} \ar@/^/[r]^{x_1} \ar@/^/[l]^{y_1}
& \ar@(ul,ur)^{z_1} \ar@/^/[r]^{x_1} \ar@/^/[l]^{y_1}
& v_{\alpha-1}  \ar@(ul,ur)^{z_1} \ar@/^/[r]^{x_1} \ar@/^/[l]^{y_1}
& v_\alpha      \ar@(ul,ur)^{z_1} \ar@/^/[r]^{x_1} \ar@/^/[l]^{y_1}
& v_{\alpha+1}  \ar@(ul,ur)^{z_1} \ar@/^/[r]^{x_1} \ar@/^/[l]^{y_1}
& v_{\alpha+2}  \ar@(ul,ur)^{z_1} \ar@/^/[r]^{x_1} \ar@/^/[l]^{y_1}
&                                    \ar@/^/[l]^{y_1} & \dots\\
}
\]
\normalsize
where we have added in the scheme the action of $z_1$ for the sake of completeness.
\item By Proposition \ref{iso-des-P-phi}, for all $\Psi\in\Z.\phi$, $P_\phi$ and $P_\Psi$ are isomorphic in 
$\B_1^{\bar{q},\Lambda}-\wMod$. 
\end{enumerate}
\end{subexample}

Recall from Remark \ref{N-phi} the submodule $N_\phi$ of $P_\phi$.

\begin{subexample} -- \label{exemple-1} \rm Suppose $\phi\in \Z.\1$. 
\begin{enumerate}
\item  First case: $\phi\in (-\N).\1$. 
In this case there exists a unique non negative integer $\alpha$ such that
$\phi=(-\alpha)\1$, or equivalently such that $\phi(z_1) = q_1^\alpha$. 
\begin{enumerate}
\item The action of $x_1$, $y_1$, $z_1$ can be pictured as follows: 
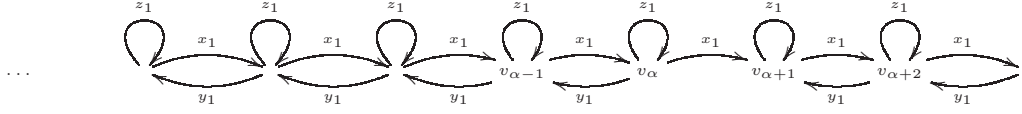
\tiny
\[
\xymatrix@!C{
\dots
& \ar@(ul,ur)^{z_1} \ar@/^/[r]^{x_1}
& \ar@(ul,ur)^{z_1} \ar@/^/[r]^{x_1} \ar@/^/[l]^{y_1}
& \ar@(ul,ur)^{z_1} \ar@/^/[r]^{x_1} \ar@/^/[l]^{y_1}
& v_{\alpha-1}  \ar@(ul,ur)^{z_1} \ar@/^/[r]^{x_1} \ar@/^/[l]^{y_1}
& v_\alpha      \ar@(ul,ur)^{z_1} \ar@/^/[r]^{x_1} \ar@/^/[l]^{y_1}
& v_{\alpha+1}  \ar@(ul,ur)^{z_1} \ar@/^/[r]^{x_1} 
& v_{\alpha+2}  \ar@(ul,ur)^{z_1} \ar@/^/[r]^{x_1} \ar@/^/[l]^{y_1}
&                                    \ar@/^/[l]^{y_1} & \dots\\
}
\]
\normalsize
\item Applying Proposition \ref{sous-objets-et-poids}, one gets at once 
that 
\[
N_\phi = \bigoplus_{k\ge \alpha +1} \k.v_k. 
\]
In particular the set of weights of $N_\phi$ is $\N^\ast.\1$.
\item  It follows that the set of weights of $S_\phi$ is $(-\N).\1$.
\end{enumerate}
\item Second case: $\phi\in \N^\ast.\1$. In this case there exists a unique positive integer 
$\beta$ such that $\phi=\beta.\1$, or equivalently such that $\phi(z_1) = q_1^{-\beta}$. 
\begin{enumerate}
\item The action of $x_1$, $y_1$, $z_1$ can be pictured as follows: 
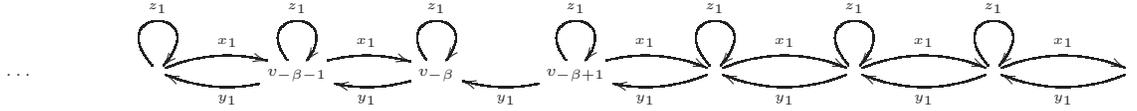
\tiny
\[
\xymatrix@!C{
\dots 
& \ar@(ul,ur)^{z_1} \ar@/^/[r]^{x_1} 
& v_{-\beta-1} \ar@(ul,ur)^{z_1} \ar@/^/[r]^{x_1} \ar@/^/[l]^{y_1}
& v_{-\beta}     \ar@(ul,ur)^{z_1}                  \ar@/^/[l]^{y_1}
& v_{-\beta+1} \ar@(ul,ur)^{z_1} \ar@/^/[r]^{x_1} \ar@/^/[l]^{y_1}
& \ar@(ul,ur)^{z_1} \ar@/^/[r]^{x_1} \ar@/^/[l]^{y_1}
& \ar@(ul,ur)^{z_1} \ar@/^/[r]^{x_1} \ar@/^/[l]^{y_1}
& \ar@(ul,ur)^{z_1} \ar@/^/[r]^{x_1} \ar@/^/[l]^{y_1}
&                                    \ar@/^/[l]^{y_1}
& \dots 
}
\]
\normalsize
\item Arguing as in the first case, $N_\phi = \bigoplus_{k\le -\beta} \k.v_k$. 
Hence, the set of weights of $N_\phi$ is $(-\N).\1$.
\item As a consequence, the set of weights of $S_\phi$ is $\N^\ast.\1$.
\end{enumerate}
\item By Proposition \ref{iso-des-P-phi}, we get that: 
\begin{enumerate}
\item for all $\rho,\psi\in(-\N).\1$, $P_\rho$ and $P_\psi$ are isomorphic in 
$\B_1^{\bar{q},\Lambda}-\wMod$;
\item for all $\rho,\psi\in\N^\ast.\1$, $P_\rho$ and $P_\psi$ are isomorphic in 
$\B_1^{\bar{q},\Lambda}-\wMod$.
\item Now, consider $\rho\in(-\N).\1$ and $\psi\in\N^\ast.\1$. 
From the analysis of the first two cases, it follows that the sets of weights 
of $N_\rho$ and $N_\psi$ are $\N^\ast.\1$ and $(-\N).\1$, respectively. 
Since we are in the generic setting, these sets are not equal. Hence, $N_\rho$ and $N_\psi$ are not
isomorphic, implying that $P_\rho$ and $P_\psi$ are not isomorphic, neither.
Using again both previous cases, the sets of weights 
of $S_\rho$ and $S_\psi$ are $(-\N).\1$ and $\N^\ast.\1$, respectively. 
Thus, $S_\rho$ and $S_\psi$ are not isomorphic.
\end{enumerate}
\end{enumerate}
\end{subexample}

The next two theorems provide a description in terms of orbits of the isomorphism classes of projective and simple weight modules, respectively. 

\begin{subtheorem} -- Let $\phi,\psi\in \widehat{R^\circ}$.
\begin{enumerate}
\item If the modules $P_\phi$ and $P_\psi$ are isomorphic, either in 
$\B_1^{\bar{q},\Lambda}-\wMod$ or in $\B_1^{\bar{q},\Lambda}-\Mod$, then
$\Z.\phi = \Z.\psi$.
\item If $\Z.\phi = \Z.\psi \neq\Z.\1$, then $P_\phi$ and $P_\psi$ are isomorphic, both in
$\B_1^{\bar{q},\Lambda}-\wMod$ and in $\B_1^{\bar{q},\Lambda}-\Mod$.
\item If $\Z.\phi = \Z.\psi = \Z.\1$, then $P_\phi$ and $P_\psi$ are isomorphic if and only if both $\phi$ and $\psi$ 
belong to $(-\N).\1$ or they both belong to $\N^\ast.\1$.
\end{enumerate}
\end{subtheorem}

\proof The proof of the first item is clear since $\Z.\phi$ and $\Z.\psi$ are the sets of weights of $P_\phi$ and $P_\psi$, respectively.
The other statements have already been proved.\qed

\begin{subtheorem} -- \label{iso-des-S}
Let $\phi,\psi\in \widehat{R^\circ}$.
\begin{enumerate}
\item  If $S_\phi$ and $S_\psi$ are isomorphic, either in $\B_1^{\bar{q},\Lambda}-\wMod$ or in $\B_1^{\bar{q},\Lambda}-\Mod$, then
$\Z.\phi = \Z.\psi$.
\item  If $\Z.\phi = \Z.\psi \neq\Z.\1$, then $S_\phi$ and $S_\psi$ are isomorphic.
\item  If $\Z.\phi = \Z.\psi = \Z.\1$, then $S_\phi$ and $S_\psi$ are isomorphic if and only if both $\phi$ and $\psi$ 
belong to $(-\N).\1$ or both of them belong to $\N^\ast.\1$.
\end{enumerate}
\end{subtheorem}

\proof Again the proof of the first item is clear since, by the previous results, 
the sets of weights of $S_\phi$ and $S_\psi$ are, respectively, subsets of $\Z.\phi$ and $\Z.\psi$. 
All the rest has already been proved.\qed

\begin{subcorollary} -- \label{iso-S-via-poids}
Let $\phi,\psi\in\widehat{R^\circ}$. The following are equivalent:
\begin{enumerate}
\item $S_\phi$ and $S_\psi$ are isomorphic either in $\B_1^{\bar{q},\Lambda}-\wMod$ or 
in $\B_1^{\bar{q},\Lambda}-\Mod$;
\item $S_\phi$ and $S_\psi$ have the same set of weights.
\end{enumerate}
\end{subcorollary}

\proof Clearly, the first statement implies the second one. 
Now, let $w(S_\phi)$ and $w(S_\psi)$ be the sets of weights of $S_\phi$ and $S_\psi$, respectively, and suppose 
$w(S_\phi)=w(S_\psi)$. There are different possible cases:\\
{\em First case:} $\phi\notin\Z\1$. Example \ref{exemple-hors-1} shows that $S_\phi=P_\phi$ and $\Z\phi=w(S_\phi)=w(S_\psi)$.
But $\psi\in w(S_\psi)$, hence $\Z\phi=\Z\psi$ and Theorem \ref{iso-des-S} proves that $S_\phi$ and $S_\psi$ are isomorphic.\\
{\em Second case:} $\phi\in(-\N)\1$. By Example \ref{exemple-1}, $(-\N)\1=w(S_\phi)=w(S_\psi)$. Hence, $\psi\in(-\N)\1$ and $\Z\phi=\Z\1$. Finally, Theorem \ref{iso-des-S} proves that $S_\phi$ and $S_\psi$ are isomorphic.\\
{\em Third case:} $\phi\in\N^\ast\1$. The same argument as in the second case applies to prove that
$S_\phi$ and $S_\psi$ are isomorphic.\qed

\subsection{Simple weight modules for arbitrary $n$}

Recall from the introduction of this section that $\B_n^{\bar{q},\Lambda}$ is identified to the $n$-fold tensor product of algebras $\B_1^{q_1} \otimes\dots\otimes \B_1^{q_n}$ and that, likewise, $R^\circ$ 
identifies to the tensor product of $n$ copies of $\k[z_1^{\pm 1}]$. In the rest of this section, 
we will freely use these identifications. 

Notice further that we are in position to apply Proposition \ref{weight-in-tpa}.

Let $\phi\in\widehat{R}$. By Proposition \ref{weight-in-tpa}, under the above identification, there is a unique $n$-tuple $(\phi_1,\dots,\phi_n)$  
of elements of $\widehat{\k[z_1^{\pm 1}]}$ such that $\phi=\phi_1\dots\phi_n$. 

We are interested in two weight $\B_n^{\bar{q},\Lambda}$-modules, namely 
$S_\phi$ and $S_{\phi_1} \otimes\dots\otimes S_{\phi_n}.$

\begin{sublemma} -- \label{tenseur-S-simple}
Let $(\phi_1,\dots,\phi_n)$ be a $n$-tuple of elements of $\widehat{\k[z_1^{\pm 1}]}$.
The $\B_n^{\bar{q},\Lambda}$-module $S_{\phi_1} \otimes\dots\otimes S_{\phi_n}$ is simple. 
\end{sublemma}

\proof Since $\bar{q}$ is generic, according to Proposition \ref{generic-weight-decomposition}, for any $i$, $1 \le i \le n$,
the weight spaces of the $\B_1^{q_i}$-module $P_{\phi_i}$ are at most 
one dimensional $\k$-subspaces; hence, the same holds for the
quotient $S_{\phi_i}$. Hence, 
the statement follows from Lemma \ref{shurian-simple}.\qed\\



\begin{subtheorem} -- \label{thm-simple}
Suppose $\Lambda$ is the skew-symmetric $n\times n$-matrix with all entries equal to $1$ and suppose $\bar{q}$ is generic.
\begin{enumerate}
\item For all $n$-tuple $(\phi_1,\dots,\phi_n)$ of elements of $\widehat{\k[z_1^{\pm 1}}]$, the  $\B_n^{\bar{q},\Lambda}$
module $S_{\phi_1} \otimes\dots\otimes S_{\phi_n}$
is a simple weight module and any simple weight $\B_n^{\bar{q},\Lambda}$-module is of that form. 
\item Let $(\phi_1,\dots,\phi_n)$ and $(\psi_1,\dots,\psi_n)$ be elements of $\left(\widehat{\k[z_1^{\pm 1}}]\right)^n$. 
The $\B_n^{\bar{q},\Lambda}$-modules $S_{\phi_1} \otimes\dots\otimes S_{\phi_n}$ and 
$S_{\psi_1} \otimes\dots\otimes S_{\psi_n}$ are isomorphic as weight modules if and only if, for all $1 \le i \le n$, 
$S_{\phi_i}$ and $S_{\psi_i}$ are isomorphic in $\B_1^{q_i} -\wMod$.
\end{enumerate}
\end{subtheorem}

\proof The proof of the first item follows from Lemma \ref{tenseur-S-simple}, Proposition \ref{ubiquite-des-P},
 and Corollary \ref{iso-L-tenseur-L}. 
We now prove the second statement. The {\em if} part of the equivalence is evident. Suppose now that 
$S_{\phi_1} \otimes\dots\otimes S_{\phi_n}$ and $S_{\psi_1} \otimes\dots\otimes S_{\psi_n}$ are isomorphic. Given $i$ such that 
$1 \le i \le n$, let $w(S_{\phi_i})$ and $w(S_{\psi_i})$ be the sets of weights of $S_{\phi_i}$ and $S_{\psi_i}$, 
respectively.
It is clear that the set of weights of $S_{\phi_1} \otimes\dots\otimes S_{\phi_n}$ is the image under the bijective map
$\left(\widehat{\k[z_1^{\pm 1}}]\right)^n \longrightarrow \widehat{R^\circ}$, 
$(\chi_1,\dots,\chi_n) \mapsto \chi_1\dots\chi_n$ --as defined in the 
second item of Proposition \ref{weight-in-tpa}-- of the subset 
$w(S_{\phi_1})\times\dots\times w(S_{\phi_1})$.
Of course, the same holds for $S_{\psi_1} \otimes\dots\otimes S_{\psi_n}$. Hence, these images must be equal
since $S_{\phi_1} \otimes\dots\otimes S_{\phi_n}$ and $S_{\psi_1} \otimes\dots\otimes S_{\psi_n}$ are isomorphic.
But that latter map is bijective. So, for $1 \le i \le n$, $w(S_{\phi_i})=w(S_{\psi_i})$. Corollary \ref{iso-S-via-poids}
yields the desired conclusion. \qed


\subsection{The case $n=2$}

In this short section we concentrate on the case $n=2$. 
Our aim is to illustrate how the case of an arbitrary $n$ can be understood on the 
basis of the case $n=1$. 
We do this by using graphs attached to the modules under consideration,
following the approach of Section \ref{exemple-n=1}. The action of $x_1$ is displayed horizontaly
from left to right, the action of $y_1$ is displayed horizontally from right to left, 
the action of $x_2$ is displayed vertically from bottom to top and the action of $y_2$ is displayed vertically from top to bottom.   

\begin{subexample} -- \rm Let $\phi \in \widehat{\k[z_1^{\pm 1},z_2^{\pm 1}]}$ and suppose 
$\comp(\phi)=\emptyset$. 

In this case, the action of $x_1$, $y_1$, $x_2$, $y_2$ can be pictured as follows:

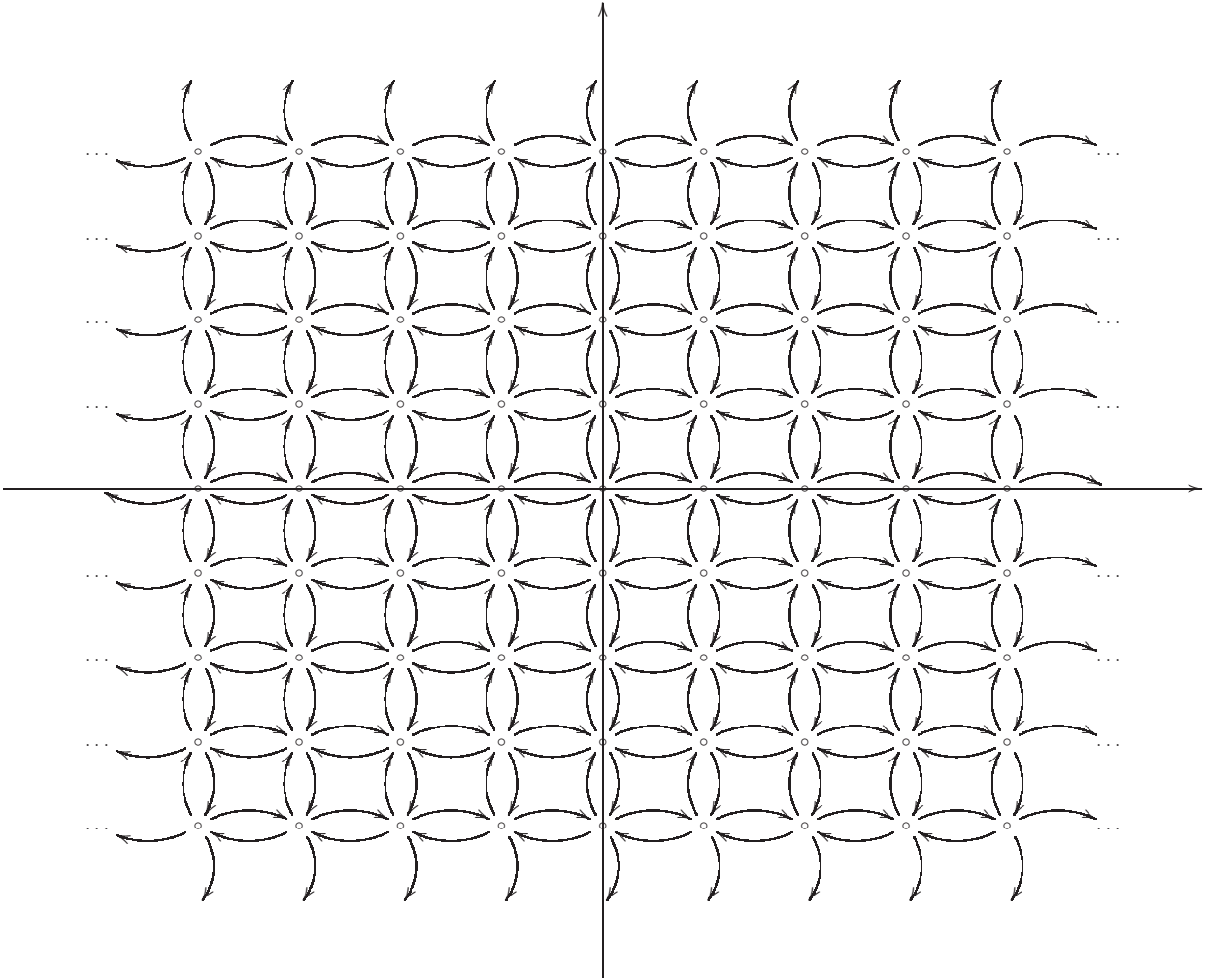
\tiny
\[
\xymatrix@!C{
&&&&&&&&&&& \cr
&&&&&&&&&&& \cr
& \dots 
& \circ \ar@/^/[r]^{} \ar@/^/[l]^{} \ar@/^/[u]^{} \ar@/^/[d]^{}  
& \circ \ar@/^/[r]^{} \ar@/^/[l]^{} \ar@/^/[u]^{} \ar@/^/[d]^{}  
& \circ \ar@/^/[r]^{} \ar@/^/[l]^{} \ar@/^/[u]^{} \ar@/^/[d]^{}  
& \circ \ar@/^/[r]^{} \ar@/^/[l]^{} \ar@/^/[u]^{} \ar@/^/[d]^{}  
& \circ \ar@/^/[r]^{} \ar@/^/[l]^{} \ar@/^/[u]^{} \ar@/^/[d]^{}  
& \circ \ar@/^/[r]^{} \ar@/^/[l]^{} \ar@/^/[u]^{} \ar@/^/[d]^{}  
& \circ \ar@/^/[r]^{} \ar@/^/[l]^{} \ar@/^/[u]^{} \ar@/^/[d]^{}  
& \circ \ar@/^/[r]^{} \ar@/^/[l]^{} \ar@/^/[u]^{} \ar@/^/[d]^{}  
& \circ \ar@/^/[r]^{} \ar@/^/[l]^{} \ar@/^/[u]^{} \ar@/^/[d]^{}  
& \dots & \\
& \dots 
& \circ \ar@/^/[r]^{} \ar@/^/[l]^{} \ar@/^/[u]^{} \ar@/^/[d]^{}  
& \circ \ar@/^/[r]^{} \ar@/^/[l]^{} \ar@/^/[u]^{} \ar@/^/[d]^{}  
& \circ \ar@/^/[r]^{} \ar@/^/[l]^{} \ar@/^/[u]^{} \ar@/^/[d]^{}  
& \circ \ar@/^/[r]^{} \ar@/^/[l]^{} \ar@/^/[u]^{} \ar@/^/[d]^{}  
& \circ \ar@/^/[r]^{} \ar@/^/[l]^{} \ar@/^/[u]^{} \ar@/^/[d]^{}  
& \circ \ar@/^/[r]^{} \ar@/^/[l]^{} \ar@/^/[u]^{} \ar@/^/[d]^{}  
& \circ \ar@/^/[r]^{} \ar@/^/[l]^{} \ar@/^/[u]^{} \ar@/^/[d]^{}  
& \circ \ar@/^/[r]^{} \ar@/^/[l]^{} \ar@/^/[u]^{} \ar@/^/[d]^{}  
& \circ \ar@/^/[r]^{} \ar@/^/[l]^{} \ar@/^/[u]^{} \ar@/^/[d]^{}  
& \dots & \\
& \dots 
& \circ \ar@/^/[r]^{} \ar@/^/[l]^{} \ar@/^/[u]^{} \ar@/^/[d]^{}
& \circ \ar@/^/[r]^{} \ar@/^/[l]^{} \ar@/^/[u]^{} \ar@/^/[d]^{} 
& \circ \ar@/^/[r]^{} \ar@/^/[l]^{} \ar@/^/[u]^{} \ar@/^/[d]^{}
& \circ \ar@/^/[r]^{} \ar@/^/[l]^{} \ar@/^/[u]^{} \ar@/^/[d]^{}  
& \circ \ar@/^/[r]^{} \ar@/^/[l]^{} \ar@/^/[u]^{} \ar@/^/[d]^{}  
& \circ \ar@/^/[r]^{} \ar@/^/[l]^{} \ar@/^/[u]^{} \ar@/^/[d]^{}  
& \circ \ar@/^/[r]^{} \ar@/^/[l]^{} \ar@/^/[u]^{} \ar@/^/[d]^{}  
& \circ \ar@/^/[r]^{} \ar@/^/[l]^{} \ar@/^/[u]^{} \ar@/^/[d]^{}  
& \circ \ar@/^/[r]^{} \ar@/^/[l]^{} \ar@/^/[u]^{} \ar@/^/[d]^{}  
& \dots & \\
& \dots 
& \circ \ar@/^/[r]^{} \ar@/^/[l]^{} \ar@/^/[u]^{} \ar@/^/[d]^{}  
& \circ \ar@/^/[r]^{} \ar@/^/[l]^{} \ar@/^/[u]^{} \ar@/^/[d]^{}  
& \circ \ar@/^/[r]^{} \ar@/^/[l]^{} \ar@/^/[u]^{} \ar@/^/[d]^{}  
& \circ \ar@/^/[r]^{} \ar@/^/[l]^{} \ar@/^/[u]^{} \ar@/^/[d]^{}  
& \circ \ar@/^/[r]^{} \ar@/^/[l]^{} \ar@/^/[u]^{} \ar@/^/[d]^{}  
& \circ \ar@/^/[r]^{} \ar@/^/[l]^{} \ar@/^/[u]^{} \ar@/^/[d]^{}  
& \circ \ar@/^/[r]^{} \ar@/^/[l]^{} \ar@/^/[u]^{} \ar@/^/[d]^{}  
& \circ \ar@/^/[r]^{} \ar@/^/[l]^{} \ar@/^/[u]^{} \ar@/^/[d]^{}  
& \circ \ar@/^/[r]^{} \ar@/^/[l]^{} \ar@/^/[u]^{} \ar@/^/[d]^{}  
& \dots & \\
\ar[rrrrrrrrrrrr]
& 
& \circ \ar@/^/[r]^{} \ar@/^/[l]^{} \ar@/^/[u]^{} \ar@/^/[d]^{}  
& \circ \ar@/^/[r]^{} \ar@/^/[l]^{} \ar@/^/[u]^{} \ar@/^/[d]^{}  
& \circ \ar@/^/[r]^{} \ar@/^/[l]^{} \ar@/^/[u]^{} \ar@/^/[d]^{}  
& \circ \ar@/^/[r]^{} \ar@/^/[l]^{} \ar@/^/[u]^{} \ar@/^/[d]^{}  
& \circ \ar@/^/[r]^{} \ar@/^/[l]^{} \ar@/^/[u]^{} \ar@/^/[d]^{}  
& \circ \ar@/^/[r]^{} \ar@/^/[l]^{} \ar@/^/[u]^{} \ar@/^/[d]^{}  
& \circ \ar@/^/[r]^{} \ar@/^/[l]^{} \ar@/^/[u]^{} \ar@/^/[d]^{}  
& \circ \ar@/^/[r]^{} \ar@/^/[l]^{} \ar@/^/[u]^{} \ar@/^/[d]^{}  
& \circ \ar@/^/[r]^{} \ar@/^/[l]^{} \ar@/^/[u]^{} \ar@/^/[d]^{}  
&  & \\
& \dots 
& \circ \ar@/^/[r]^{} \ar@/^/[l]^{} \ar@/^/[u]^{} \ar@/^/[d]^{}  
& \circ \ar@/^/[r]^{} \ar@/^/[l]^{} \ar@/^/[u]^{} \ar@/^/[d]^{}  
& \circ \ar@/^/[r]^{} \ar@/^/[l]^{} \ar@/^/[u]^{} \ar@/^/[d]^{}  
& \circ \ar@/^/[r]^{} \ar@/^/[l]^{} \ar@/^/[u]^{} \ar@/^/[d]^{}  
& \circ \ar@/^/[r]^{} \ar@/^/[l]^{} \ar@/^/[u]^{} \ar@/^/[d]^{}  
& \circ \ar@/^/[r]^{} \ar@/^/[l]^{} \ar@/^/[u]^{} \ar@/^/[d]^{}  
& \circ \ar@/^/[r]^{} \ar@/^/[l]^{} \ar@/^/[u]^{} \ar@/^/[d]^{}  
& \circ \ar@/^/[r]^{} \ar@/^/[l]^{} \ar@/^/[u]^{} \ar@/^/[d]^{}  
& \circ \ar@/^/[r]^{} \ar@/^/[l]^{} \ar@/^/[u]^{} \ar@/^/[d]^{}  
& \dots & \\
& \dots 
& \circ \ar@/^/[r]^{} \ar@/^/[l]^{} \ar@/^/[u]^{} \ar@/^/[d]^{}  
& \circ \ar@/^/[r]^{} \ar@/^/[l]^{} \ar@/^/[u]^{} \ar@/^/[d]^{}  
& \circ \ar@/^/[r]^{} \ar@/^/[l]^{} \ar@/^/[u]^{} \ar@/^/[d]^{}  
& \circ \ar@/^/[r]^{} \ar@/^/[l]^{} \ar@/^/[u]^{} \ar@/^/[d]^{}  
& \circ \ar@/^/[r]^{} \ar@/^/[l]^{} \ar@/^/[u]^{} \ar@/^/[d]^{}  
& \circ \ar@/^/[r]^{} \ar@/^/[l]^{} \ar@/^/[u]^{} \ar@/^/[d]^{}  
& \circ \ar@/^/[r]^{} \ar@/^/[l]^{} \ar@/^/[u]^{} \ar@/^/[d]^{}  
& \circ \ar@/^/[r]^{} \ar@/^/[l]^{} \ar@/^/[u]^{} \ar@/^/[d]^{}  
& \circ \ar@/^/[r]^{} \ar@/^/[l]^{} \ar@/^/[u]^{} \ar@/^/[d]^{}  
& \dots & \\
& \dots 
& \circ \ar@/^/[r]^{} \ar@/^/[l]^{} \ar@/^/[u]^{} \ar@/^/[d]^{}  
& \circ \ar@/^/[r]^{} \ar@/^/[l]^{} \ar@/^/[u]^{} \ar@/^/[d]^{}  
& \circ \ar@/^/[r]^{} \ar@/^/[l]^{} \ar@/^/[u]^{} \ar@/^/[d]^{}  
& \circ \ar@/^/[r]^{} \ar@/^/[l]^{} \ar@/^/[u]^{} \ar@/^/[d]^{}  
& \circ \ar@/^/[r]^{} \ar@/^/[l]^{} \ar@/^/[u]^{} \ar@/^/[d]^{}  
& \circ \ar@/^/[r]^{} \ar@/^/[l]^{} \ar@/^/[u]^{} \ar@/^/[d]^{}  
& \circ \ar@/^/[r]^{} \ar@/^/[l]^{} \ar@/^/[u]^{} \ar@/^/[d]^{}  
& \circ \ar@/^/[r]^{} \ar@/^/[l]^{} \ar@/^/[u]^{} \ar@/^/[d]^{}  
& \circ \ar@/^/[r]^{} \ar@/^/[l]^{} \ar@/^/[u]^{} \ar@/^/[d]^{}  
& \dots & \\
& \dots 
& \circ \ar@/^/[r]^{} \ar@/^/[l]^{} \ar@/^/[u]^{} \ar@/^/[d]^{}  
& \circ \ar@/^/[r]^{} \ar@/^/[l]^{} \ar@/^/[u]^{} \ar@/^/[d]^{}  
& \circ \ar@/^/[r]^{} \ar@/^/[l]^{} \ar@/^/[u]^{} \ar@/^/[d]^{}  
& \circ \ar@/^/[r]^{} \ar@/^/[l]^{} \ar@/^/[u]^{} \ar@/^/[d]^{}  
& \circ \ar@/^/[r]^{} \ar@/^/[l]^{} \ar@/^/[u]^{} \ar@/^/[d]^{}  
& \circ \ar@/^/[r]^{} \ar@/^/[l]^{} \ar@/^/[u]^{} \ar@/^/[d]^{}  
& \circ \ar@/^/[r]^{} \ar@/^/[l]^{} \ar@/^/[u]^{} \ar@/^/[d]^{}  
& \circ \ar@/^/[r]^{} \ar@/^/[l]^{} \ar@/^/[u]^{} \ar@/^/[d]^{}  
& \circ \ar@/^/[r]^{} \ar@/^/[l]^{} \ar@/^/[u]^{} \ar@/^/[d]^{}  
& \dots & \\
&&&&&&&&&&& \cr
&&&&&& \ar[uuuuuuuuuuuu]&&&&&}
\]
\normalsize
It is clear that in this graph, there is a path going from any vertex to the vertex $v_0$. 
It follows that any non-trivial submodule of $P_\phi$ must contain $v_0$ and hence must 
coincide with $P_\phi$. So, $P_\phi$ is simple.
\end{subexample}

\newpage

\begin{subexample} -- \rm Let $\phi \in \widehat{\k[z_1^{\pm 1},z_2^{\pm 1}]}$ and suppose 
$\comp(\phi)=\{1\}$. In this case, there exists an integer $\alpha_1$ such that
$\phi(z_1) = q_1^{\alpha_1}$. 

Let us suppose that $\alpha_1 \in \N$. 
The action of $x_1$, $y_1$, $x_2$, $y_2$ can be pictured as follows: \\

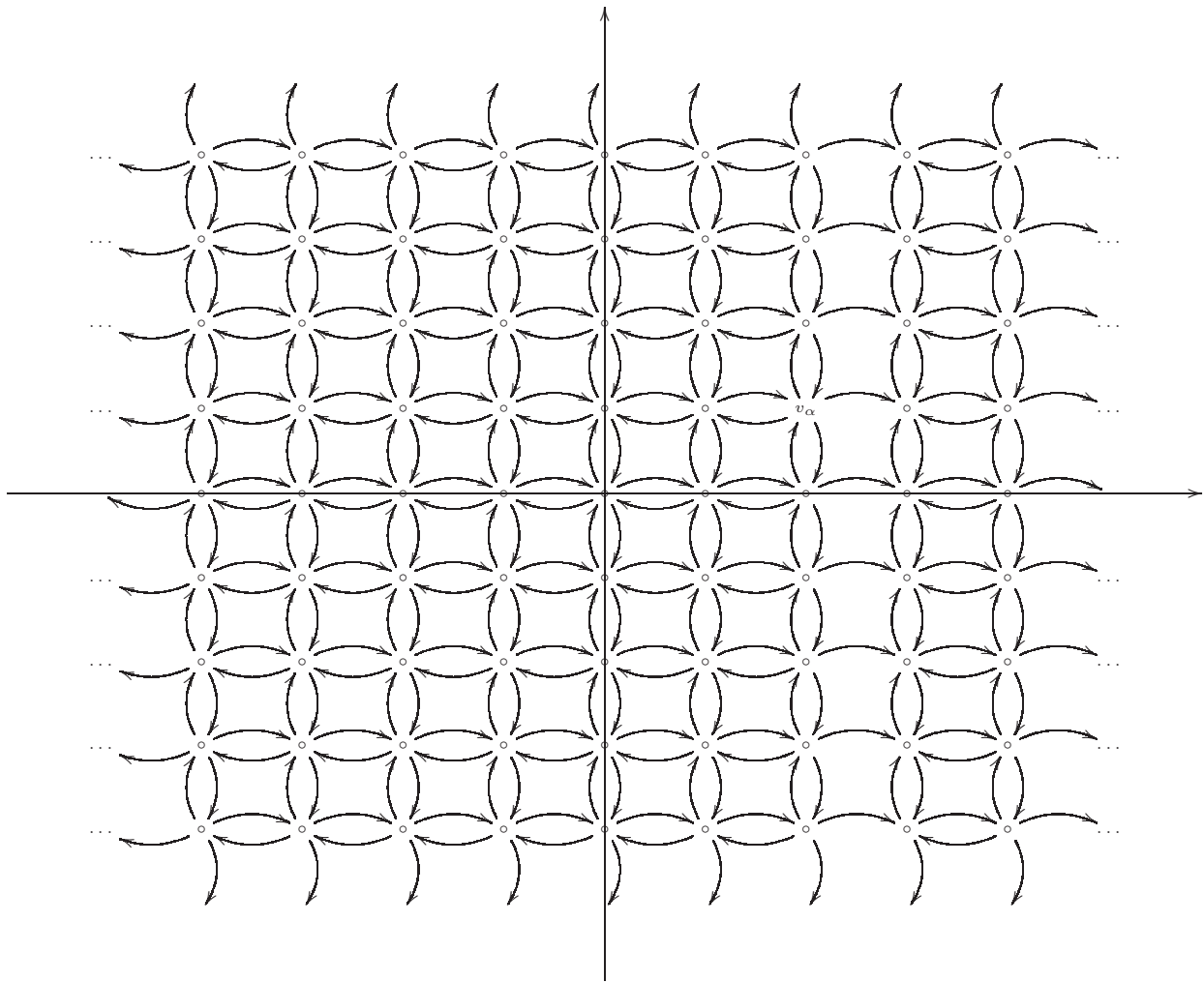
\tiny
\[
\xymatrix@!C{
&&&&&&&&&&& \cr
&&&&&&&&&&& \cr
& \dots 
& \circ \ar@/^/[r]^{} \ar@/^/[l]^{} \ar@/^/[u]^{} \ar@/^/[d]^{}  
& \circ \ar@/^/[r]^{} \ar@/^/[l]^{} \ar@/^/[u]^{} \ar@/^/[d]^{}  
& \circ \ar@/^/[r]^{} \ar@/^/[l]^{} \ar@/^/[u]^{} \ar@/^/[d]^{}  
& \circ \ar@/^/[r]^{} \ar@/^/[l]^{} \ar@/^/[u]^{} \ar@/^/[d]^{}  
& \circ \ar@/^/[r]^{} \ar@/^/[l]^{} \ar@/^/[u]^{} \ar@/^/[d]^{}  
& \circ \ar@/^/[r]^{} \ar@/^/[l]^{} \ar@/^/[u]^{} \ar@/^/[d]^{}  
& \circ \ar@/^/[r]^{} \ar@/^/[l]^{} \ar@/^/[u]^{} \ar@/^/[d]^{}  
& \circ \ar@/^/[r]^{}               \ar@/^/[u]^{} \ar@/^/[d]^{}  
& \circ \ar@/^/[r]^{} \ar@/^/[l]^{} \ar@/^/[u]^{} \ar@/^/[d]^{}  
& \dots & \\
& \dots 
& \circ \ar@/^/[r]^{} \ar@/^/[l]^{} \ar@/^/[u]^{} \ar@/^/[d]^{}  
& \circ \ar@/^/[r]^{} \ar@/^/[l]^{} \ar@/^/[u]^{} \ar@/^/[d]^{}  
& \circ \ar@/^/[r]^{} \ar@/^/[l]^{} \ar@/^/[u]^{} \ar@/^/[d]^{}  
& \circ \ar@/^/[r]^{} \ar@/^/[l]^{} \ar@/^/[u]^{} \ar@/^/[d]^{}  
& \circ \ar@/^/[r]^{} \ar@/^/[l]^{} \ar@/^/[u]^{} \ar@/^/[d]^{}  
& \circ \ar@/^/[r]^{} \ar@/^/[l]^{} \ar@/^/[u]^{} \ar@/^/[d]^{}  
& \circ \ar@/^/[r]^{} \ar@/^/[l]^{} \ar@/^/[u]^{} \ar@/^/[d]^{}  
& \circ \ar@/^/[r]^{}               \ar@/^/[u]^{} \ar@/^/[d]^{}  
& \circ \ar@/^/[r]^{} \ar@/^/[l]^{} \ar@/^/[u]^{} \ar@/^/[d]^{}  
& \dots & \\
& \dots 
& \circ \ar@/^/[r]^{} \ar@/^/[l]^{} \ar@/^/[u]^{} \ar@/^/[d]^{}
& \circ \ar@/^/[r]^{} \ar@/^/[l]^{} \ar@/^/[u]^{} \ar@/^/[d]^{} 
& \circ \ar@/^/[r]^{} \ar@/^/[l]^{} \ar@/^/[u]^{} \ar@/^/[d]^{}
& \circ \ar@/^/[r]^{} \ar@/^/[l]^{} \ar@/^/[u]^{} \ar@/^/[d]^{}  
& \circ \ar@/^/[r]^{} \ar@/^/[l]^{} \ar@/^/[u]^{} \ar@/^/[d]^{}  
& \circ \ar@/^/[r]^{} \ar@/^/[l]^{} \ar@/^/[u]^{} \ar@/^/[d]^{}  
& \circ \ar@/^/[r]^{} \ar@/^/[l]^{} \ar@/^/[u]^{} \ar@/^/[d]^{}  
& \circ \ar@/^/[r]^{}               \ar@/^/[u]^{} \ar@/^/[d]^{}  
& \circ \ar@/^/[r]^{} \ar@/^/[l]^{} \ar@/^/[u]^{} \ar@/^/[d]^{}  
& \dots & \\
& \dots 
& \circ    \ar@/^/[r]^{} \ar@/^/[l]^{} \ar@/^/[u]^{} \ar@/^/[d]^{}  
& \circ    \ar@/^/[r]^{} \ar@/^/[l]^{} \ar@/^/[u]^{} \ar@/^/[d]^{}  
& \circ    \ar@/^/[r]^{} \ar@/^/[l]^{} \ar@/^/[u]^{} \ar@/^/[d]^{}  
& \circ    \ar@/^/[r]^{} \ar@/^/[l]^{} \ar@/^/[u]^{} \ar@/^/[d]^{}  
& \circ    \ar@/^/[r]^{} \ar@/^/[l]^{} \ar@/^/[u]^{} \ar@/^/[d]^{}  
& \circ    \ar@/^/[r]^{} \ar@/^/[l]^{} \ar@/^/[u]^{} \ar@/^/[d]^{}  
& v_\alpha \ar@/^/[r]^{} \ar@/^/[l]^{} \ar@/^/[u]^{} \ar@/^/[d]^{}  
& \circ    \ar@/^/[r]^{}               \ar@/^/[u]^{} \ar@/^/[d]^{}  
& \circ    \ar@/^/[r]^{} \ar@/^/[l]^{} \ar@/^/[u]^{} \ar@/^/[d]^{}  
& \dots & \\
\ar[rrrrrrrrrrrr]
& 
& \circ \ar@/^/[r]^{} \ar@/^/[l]^{} \ar@/^/[u]^{} \ar@/^/[d]^{}  
& \circ \ar@/^/[r]^{} \ar@/^/[l]^{} \ar@/^/[u]^{} \ar@/^/[d]^{}  
& \circ \ar@/^/[r]^{} \ar@/^/[l]^{} \ar@/^/[u]^{} \ar@/^/[d]^{}  
& \circ \ar@/^/[r]^{} \ar@/^/[l]^{} \ar@/^/[u]^{} \ar@/^/[d]^{}  
& \circ \ar@/^/[r]^{} \ar@/^/[l]^{} \ar@/^/[u]^{} \ar@/^/[d]^{}  
& \circ \ar@/^/[r]^{} \ar@/^/[l]^{} \ar@/^/[u]^{} \ar@/^/[d]^{}  
& \circ \ar@/^/[r]^{} \ar@/^/[l]^{} \ar@/^/[u]^{} \ar@/^/[d]^{}  
& \circ \ar@/^/[r]^{}               \ar@/^/[u]^{} \ar@/^/[d]^{}  
& \circ \ar@/^/[r]^{} \ar@/^/[l]^{} \ar@/^/[u]^{} \ar@/^/[d]^{}  
&  & \\
& \dots 
& \circ \ar@/^/[r]^{} \ar@/^/[l]^{} \ar@/^/[u]^{} \ar@/^/[d]^{}  
& \circ \ar@/^/[r]^{} \ar@/^/[l]^{} \ar@/^/[u]^{} \ar@/^/[d]^{}  
& \circ \ar@/^/[r]^{} \ar@/^/[l]^{} \ar@/^/[u]^{} \ar@/^/[d]^{}  
& \circ \ar@/^/[r]^{} \ar@/^/[l]^{} \ar@/^/[u]^{} \ar@/^/[d]^{}  
& \circ \ar@/^/[r]^{} \ar@/^/[l]^{} \ar@/^/[u]^{} \ar@/^/[d]^{}  
& \circ \ar@/^/[r]^{} \ar@/^/[l]^{} \ar@/^/[u]^{} \ar@/^/[d]^{}  
& \circ \ar@/^/[r]^{} \ar@/^/[l]^{} \ar@/^/[u]^{} \ar@/^/[d]^{}  
& \circ \ar@/^/[r]^{}               \ar@/^/[u]^{} \ar@/^/[d]^{}  
& \circ \ar@/^/[r]^{} \ar@/^/[l]^{} \ar@/^/[u]^{} \ar@/^/[d]^{}  
& \dots & \\
& \dots 
& \circ \ar@/^/[r]^{} \ar@/^/[l]^{} \ar@/^/[u]^{} \ar@/^/[d]^{}  
& \circ \ar@/^/[r]^{} \ar@/^/[l]^{} \ar@/^/[u]^{} \ar@/^/[d]^{}  
& \circ \ar@/^/[r]^{} \ar@/^/[l]^{} \ar@/^/[u]^{} \ar@/^/[d]^{}  
& \circ \ar@/^/[r]^{} \ar@/^/[l]^{} \ar@/^/[u]^{} \ar@/^/[d]^{}  
& \circ \ar@/^/[r]^{} \ar@/^/[l]^{} \ar@/^/[u]^{} \ar@/^/[d]^{}  
& \circ \ar@/^/[r]^{} \ar@/^/[l]^{} \ar@/^/[u]^{} \ar@/^/[d]^{}  
& \circ \ar@/^/[r]^{} \ar@/^/[l]^{} \ar@/^/[u]^{} \ar@/^/[d]^{}  
& \circ \ar@/^/[r]^{}               \ar@/^/[u]^{} \ar@/^/[d]^{}  
& \circ \ar@/^/[r]^{} \ar@/^/[l]^{} \ar@/^/[u]^{} \ar@/^/[d]^{}  
& \dots & \\
& \dots 
& \circ \ar@/^/[r]^{} \ar@/^/[l]^{} \ar@/^/[u]^{} \ar@/^/[d]^{}  
& \circ \ar@/^/[r]^{} \ar@/^/[l]^{} \ar@/^/[u]^{} \ar@/^/[d]^{}  
& \circ \ar@/^/[r]^{} \ar@/^/[l]^{} \ar@/^/[u]^{} \ar@/^/[d]^{}  
& \circ \ar@/^/[r]^{} \ar@/^/[l]^{} \ar@/^/[u]^{} \ar@/^/[d]^{}  
& \circ \ar@/^/[r]^{} \ar@/^/[l]^{} \ar@/^/[u]^{} \ar@/^/[d]^{}  
& \circ \ar@/^/[r]^{} \ar@/^/[l]^{} \ar@/^/[u]^{} \ar@/^/[d]^{}  
& \circ \ar@/^/[r]^{} \ar@/^/[l]^{} \ar@/^/[u]^{} \ar@/^/[d]^{}  
& \circ \ar@/^/[r]^{}               \ar@/^/[u]^{} \ar@/^/[d]^{}  
& \circ \ar@/^/[r]^{} \ar@/^/[l]^{} \ar@/^/[u]^{} \ar@/^/[d]^{}  
& \dots & \\
& \dots 
& \circ \ar@/^/[r]^{} \ar@/^/[l]^{} \ar@/^/[u]^{} \ar@/^/[d]^{}  
& \circ \ar@/^/[r]^{} \ar@/^/[l]^{} \ar@/^/[u]^{} \ar@/^/[d]^{}  
& \circ \ar@/^/[r]^{} \ar@/^/[l]^{} \ar@/^/[u]^{} \ar@/^/[d]^{}  
& \circ \ar@/^/[r]^{} \ar@/^/[l]^{} \ar@/^/[u]^{} \ar@/^/[d]^{}  
& \circ \ar@/^/[r]^{} \ar@/^/[l]^{} \ar@/^/[u]^{} \ar@/^/[d]^{}  
& \circ \ar@/^/[r]^{} \ar@/^/[l]^{} \ar@/^/[u]^{} \ar@/^/[d]^{}  
& \circ \ar@/^/[r]^{} \ar@/^/[l]^{} \ar@/^/[u]^{} \ar@/^/[d]^{}  
& \circ \ar@/^/[r]^{}               \ar@/^/[u]^{} \ar@/^/[d]^{}  
& \circ \ar@/^/[r]^{} \ar@/^/[l]^{} \ar@/^/[u]^{} \ar@/^/[d]^{}  
& \dots & \\
&&&&&&&&&&& \cr
&&&&&& \ar[uuuuuuuuuuuu]&&&&&}
\]
\normalsize

The graph above is clearly divided into two parts, delimited by the vertical line intersecting
the $x$-axis at $x=\alpha_1 +1/2$. From any vertex located at the left hand side of this line, 
there is a path to $v_0$. Hence, any submodule containing one such vertex must coincide with
$P_\phi$. In addition, the $\k$-vector space with basis the vertices located to the right of this 
line is a proper submodule of $P_\phi$. It follows that this subspace must be maximum among non 
trivial submodules. That is, we have: 
\[
N_\phi = \bigoplus_{\{ k| k_1> \alpha_1\} } \k.v_k .
\]

We let the interested reader deal with the case where $\alpha_1 < 0$. Clearly, the case where 
$\comp(\phi)=\{2\}$ is similar, interchanging horizontal and vertical directions.
\end{subexample}

\newpage

\begin{subexample} -- \rm Let $\phi \in \widehat{\k[z_1^{\pm 1},z_2^{\pm 1}]}$ and suppose 
$\comp(\phi)=\{1,2\}$. In this case, there exist integers $\alpha_1$ and $\alpha_2$ such that
$\phi(z_1) = q_1^{\alpha_1}$, and $\phi(z_2) = q_2^{\alpha_2}$. 

Let us suppose that $\alpha_1, \alpha_2 \in \N$. 
Then, the action of $x_1$, $y_1$, $x_2$, $y_2$ can be pictured as follows, where 
$\alpha=(\alpha_1,\alpha_2)$: \\

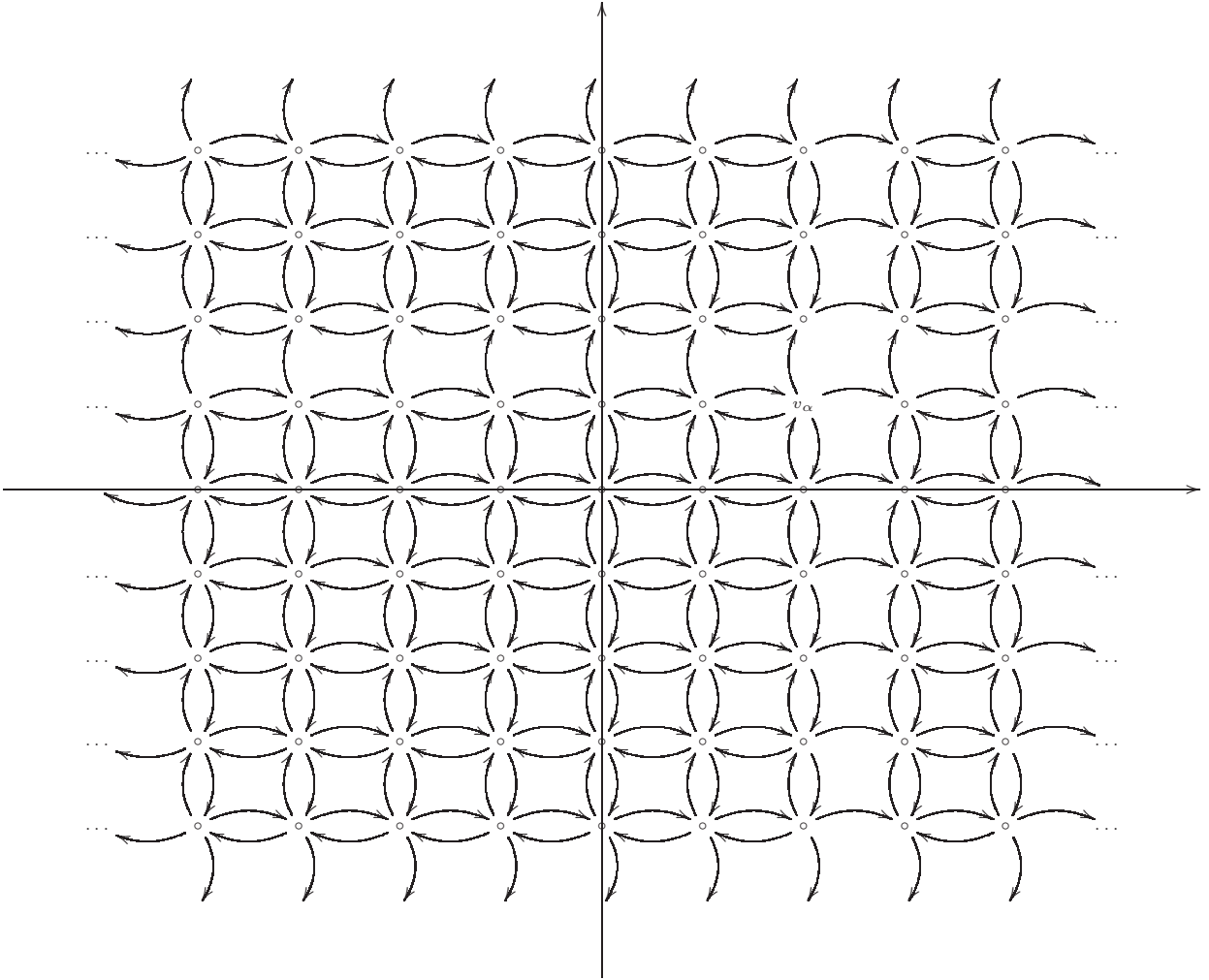
\tiny
\[
\xymatrix@!C{
&&&&&&&&&&& \cr
&&&&&&&&&&& \cr
& \dots 
& \circ \ar@/^/[r]^{} \ar@/^/[l]^{} \ar@/^/[u]^{} \ar@/^/[d]^{}  
& \circ \ar@/^/[r]^{} \ar@/^/[l]^{} \ar@/^/[u]^{} \ar@/^/[d]^{}  
& \circ \ar@/^/[r]^{} \ar@/^/[l]^{} \ar@/^/[u]^{} \ar@/^/[d]^{}  
& \circ \ar@/^/[r]^{} \ar@/^/[l]^{} \ar@/^/[u]^{} \ar@/^/[d]^{}  
& \circ \ar@/^/[r]^{} \ar@/^/[l]^{} \ar@/^/[u]^{} \ar@/^/[d]^{}  
& \circ \ar@/^/[r]^{} \ar@/^/[l]^{} \ar@/^/[u]^{} \ar@/^/[d]^{}  
& \circ \ar@/^/[r]^{} \ar@/^/[l]^{} \ar@/^/[u]^{} \ar@/^/[d]^{}  
& \circ \ar@/^/[r]^{}               \ar@/^/[u]^{} \ar@/^/[d]^{}  
& \circ \ar@/^/[r]^{} \ar@/^/[l]^{} \ar@/^/[u]^{} \ar@/^/[d]^{}  
& \dots & \\
& \dots 
& \circ \ar@/^/[r]^{} \ar@/^/[l]^{} \ar@/^/[u]^{} \ar@/^/[d]^{}  
& \circ \ar@/^/[r]^{} \ar@/^/[l]^{} \ar@/^/[u]^{} \ar@/^/[d]^{}  
& \circ \ar@/^/[r]^{} \ar@/^/[l]^{} \ar@/^/[u]^{} \ar@/^/[d]^{}  
& \circ \ar@/^/[r]^{} \ar@/^/[l]^{} \ar@/^/[u]^{} \ar@/^/[d]^{}  
& \circ \ar@/^/[r]^{} \ar@/^/[l]^{} \ar@/^/[u]^{} \ar@/^/[d]^{}  
& \circ \ar@/^/[r]^{} \ar@/^/[l]^{} \ar@/^/[u]^{} \ar@/^/[d]^{}  
& \circ \ar@/^/[r]^{} \ar@/^/[l]^{} \ar@/^/[u]^{} \ar@/^/[d]^{}  
& \circ \ar@/^/[r]^{}               \ar@/^/[u]^{} \ar@/^/[d]^{}  
& \circ \ar@/^/[r]^{} \ar@/^/[l]^{} \ar@/^/[u]^{} \ar@/^/[d]^{}  
& \dots & \\
& \dots 
& \circ \ar@/^/[r]^{} \ar@/^/[l]^{} \ar@/^/[u]^{} 
& \circ \ar@/^/[r]^{} \ar@/^/[l]^{} \ar@/^/[u]^{}  
& \circ \ar@/^/[r]^{} \ar@/^/[l]^{} \ar@/^/[u]^{} 
& \circ \ar@/^/[r]^{} \ar@/^/[l]^{} \ar@/^/[u]^{}   
& \circ \ar@/^/[r]^{} \ar@/^/[l]^{} \ar@/^/[u]^{}   
& \circ \ar@/^/[r]^{} \ar@/^/[l]^{} \ar@/^/[u]^{}   
& \circ \ar@/^/[r]^{} \ar@/^/[l]^{} \ar@/^/[u]^{}   
& \circ \ar@/^/[r]^{}               \ar@/^/[u]^{}   
& \circ \ar@/^/[r]^{} \ar@/^/[l]^{} \ar@/^/[u]^{}   
& \dots & \\
& \dots 
& \circ    \ar@/^/[r]^{} \ar@/^/[l]^{} \ar@/^/[u]^{} \ar@/^/[d]^{}  
& \circ    \ar@/^/[r]^{} \ar@/^/[l]^{} \ar@/^/[u]^{} \ar@/^/[d]^{}  
& \circ    \ar@/^/[r]^{} \ar@/^/[l]^{} \ar@/^/[u]^{} \ar@/^/[d]^{}  
& \circ    \ar@/^/[r]^{} \ar@/^/[l]^{} \ar@/^/[u]^{} \ar@/^/[d]^{}  
& \circ    \ar@/^/[r]^{} \ar@/^/[l]^{} \ar@/^/[u]^{} \ar@/^/[d]^{}  
& \circ    \ar@/^/[r]^{} \ar@/^/[l]^{} \ar@/^/[u]^{} \ar@/^/[d]^{}  
& v_\alpha \ar@/^/[r]^{} \ar@/^/[l]^{} \ar@/^/[u]^{} \ar@/^/[d]^{}  
& \circ    \ar@/^/[r]^{}               \ar@/^/[u]^{} \ar@/^/[d]^{}  
& \circ    \ar@/^/[r]^{} \ar@/^/[l]^{} \ar@/^/[u]^{} \ar@/^/[d]^{}  
& \dots & \\
\ar[rrrrrrrrrrrr]
& 
& \circ \ar@/^/[r]^{} \ar@/^/[l]^{} \ar@/^/[u]^{} \ar@/^/[d]^{}  
& \circ \ar@/^/[r]^{} \ar@/^/[l]^{} \ar@/^/[u]^{} \ar@/^/[d]^{}  
& \circ \ar@/^/[r]^{} \ar@/^/[l]^{} \ar@/^/[u]^{} \ar@/^/[d]^{}  
& \circ \ar@/^/[r]^{} \ar@/^/[l]^{} \ar@/^/[u]^{} \ar@/^/[d]^{}  
& \circ \ar@/^/[r]^{} \ar@/^/[l]^{} \ar@/^/[u]^{} \ar@/^/[d]^{}  
& \circ \ar@/^/[r]^{} \ar@/^/[l]^{} \ar@/^/[u]^{} \ar@/^/[d]^{}  
& \circ \ar@/^/[r]^{} \ar@/^/[l]^{} \ar@/^/[u]^{} \ar@/^/[d]^{}  
& \circ \ar@/^/[r]^{}               \ar@/^/[u]^{} \ar@/^/[d]^{}  
& \circ \ar@/^/[r]^{} \ar@/^/[l]^{} \ar@/^/[u]^{} \ar@/^/[d]^{}  
&  & \\
& \dots 
& \circ \ar@/^/[r]^{} \ar@/^/[l]^{} \ar@/^/[u]^{} \ar@/^/[d]^{}  
& \circ \ar@/^/[r]^{} \ar@/^/[l]^{} \ar@/^/[u]^{} \ar@/^/[d]^{}  
& \circ \ar@/^/[r]^{} \ar@/^/[l]^{} \ar@/^/[u]^{} \ar@/^/[d]^{}  
& \circ \ar@/^/[r]^{} \ar@/^/[l]^{} \ar@/^/[u]^{} \ar@/^/[d]^{}  
& \circ \ar@/^/[r]^{} \ar@/^/[l]^{} \ar@/^/[u]^{} \ar@/^/[d]^{}  
& \circ \ar@/^/[r]^{} \ar@/^/[l]^{} \ar@/^/[u]^{} \ar@/^/[d]^{}  
& \circ \ar@/^/[r]^{} \ar@/^/[l]^{} \ar@/^/[u]^{} \ar@/^/[d]^{}  
& \circ \ar@/^/[r]^{}               \ar@/^/[u]^{} \ar@/^/[d]^{}  
& \circ \ar@/^/[r]^{} \ar@/^/[l]^{} \ar@/^/[u]^{} \ar@/^/[d]^{}  
& \dots & \\
& \dots 
& \circ \ar@/^/[r]^{} \ar@/^/[l]^{} \ar@/^/[u]^{} \ar@/^/[d]^{}  
& \circ \ar@/^/[r]^{} \ar@/^/[l]^{} \ar@/^/[u]^{} \ar@/^/[d]^{}  
& \circ \ar@/^/[r]^{} \ar@/^/[l]^{} \ar@/^/[u]^{} \ar@/^/[d]^{}  
& \circ \ar@/^/[r]^{} \ar@/^/[l]^{} \ar@/^/[u]^{} \ar@/^/[d]^{}  
& \circ \ar@/^/[r]^{} \ar@/^/[l]^{} \ar@/^/[u]^{} \ar@/^/[d]^{}  
& \circ \ar@/^/[r]^{} \ar@/^/[l]^{} \ar@/^/[u]^{} \ar@/^/[d]^{}  
& \circ \ar@/^/[r]^{} \ar@/^/[l]^{} \ar@/^/[u]^{} \ar@/^/[d]^{}  
& \circ \ar@/^/[r]^{}               \ar@/^/[u]^{} \ar@/^/[d]^{}  
& \circ \ar@/^/[r]^{} \ar@/^/[l]^{} \ar@/^/[u]^{} \ar@/^/[d]^{}  
& \dots & \\
& \dots 
& \circ \ar@/^/[r]^{} \ar@/^/[l]^{} \ar@/^/[u]^{} \ar@/^/[d]^{}  
& \circ \ar@/^/[r]^{} \ar@/^/[l]^{} \ar@/^/[u]^{} \ar@/^/[d]^{}  
& \circ \ar@/^/[r]^{} \ar@/^/[l]^{} \ar@/^/[u]^{} \ar@/^/[d]^{}  
& \circ \ar@/^/[r]^{} \ar@/^/[l]^{} \ar@/^/[u]^{} \ar@/^/[d]^{}  
& \circ \ar@/^/[r]^{} \ar@/^/[l]^{} \ar@/^/[u]^{} \ar@/^/[d]^{}  
& \circ \ar@/^/[r]^{} \ar@/^/[l]^{} \ar@/^/[u]^{} \ar@/^/[d]^{}  
& \circ \ar@/^/[r]^{} \ar@/^/[l]^{} \ar@/^/[u]^{} \ar@/^/[d]^{}  
& \circ \ar@/^/[r]^{}               \ar@/^/[u]^{} \ar@/^/[d]^{}  
& \circ \ar@/^/[r]^{} \ar@/^/[l]^{} \ar@/^/[u]^{} \ar@/^/[d]^{}  
& \dots & \\
& \dots 
& \circ \ar@/^/[r]^{} \ar@/^/[l]^{} \ar@/^/[u]^{} \ar@/^/[d]^{}  
& \circ \ar@/^/[r]^{} \ar@/^/[l]^{} \ar@/^/[u]^{} \ar@/^/[d]^{}  
& \circ \ar@/^/[r]^{} \ar@/^/[l]^{} \ar@/^/[u]^{} \ar@/^/[d]^{}  
& \circ \ar@/^/[r]^{} \ar@/^/[l]^{} \ar@/^/[u]^{} \ar@/^/[d]^{}  
& \circ \ar@/^/[r]^{} \ar@/^/[l]^{} \ar@/^/[u]^{} \ar@/^/[d]^{}  
& \circ \ar@/^/[r]^{} \ar@/^/[l]^{} \ar@/^/[u]^{} \ar@/^/[d]^{}  
& \circ \ar@/^/[r]^{} \ar@/^/[l]^{} \ar@/^/[u]^{} \ar@/^/[d]^{}  
& \circ \ar@/^/[r]^{}               \ar@/^/[u]^{} \ar@/^/[d]^{}  
& \circ \ar@/^/[r]^{} \ar@/^/[l]^{} \ar@/^/[u]^{} \ar@/^/[d]^{}  
& \dots & \\
&&&&&&&&&&& \cr
&&&&&& \ar[uuuuuuuuuuuu]&&&&&}
\]
\normalsize

In this case, the graph is again divided into two parts, delimited by the vertical line intersecting
the $x$-axis at $x=\alpha_1 +1/2$ and the horizontal line intersecting the $y$-axis at 
$y = \alpha_2 + 1/2$.

From any vertex located at the left hand side of the vertical line and bellow the horizontal one, 
there is a path to $v_0$. Hence, any submodule containing one such vertex must coincide with
$P_\phi$. In addition, the $\k$-vector space with basis the vertices located either to the right of 
the vertical line or above the horizontal one is a submodule of $P_\phi$. 
It follows that this subspace must be maximum among non 
trivial submodules. That is, we have: 
\[
N_\phi = \bigoplus_{\{k| k_1 > \alpha_1 {\text or} k_2 > \alpha_2\} } \k.v_k .
\]
\end{subexample}

\section{Quantized Weyl algebras and Zhang twists} \label{QWa-et-ZT}

The aim of this section is to show how weight modules of the extension 
$(R^\circ,\B_n^{\bar{q},\Lambda})$ may be described using their analogues of the special case
where $\Lambda$ has all its entries equal to $1$. The method that we use to reduce the general case 
to this special case is based on the notion of Zhang twist as introduced in \cite{Z}. 

We first recall this notion, then show that $\B_n^{\bar{q},\Lambda}$ is indeed a Zhang twist of 
$\B_n^{\bar{q},(1)}$ and finally establish that all simple weight modules 
over $\B_n^{\bar{q},\Lambda}$ are in fact obtained as Zhang twists of simple weight modules over 
$\B_n^{\bar{q},(1)}$.\\

For all this section, we fix $\bar{q}=(q_1,\dots,q_n) \in (\k^\ast)^n$ and 
a skew-symmetric $n \times n$ matrix $\Lambda=(\lambda_{ij})$.
It is clear from the definition of 
$\A_n^{\bar{q},\Lambda}$ that it is $\Z^n$-graded by 
\[
\deg(x_i) = e_i
\qquad\mbox{and}\qquad
\deg(y_i) = -e_i .
\]
From Proposition \ref{R-basis-of-A}, we get that the homogeneous part of degree $0$ with respect to this grading is reduced to 
\[
R = \k[z_1,\dots,z_n]; 
\]
the subalgebra of $\A_n^{\bar{q},\Lambda}$ generated by $z_i$, $1 \le i \le n$, the latter being algebraically independent.

\subsection{Some $\k$-algebra automorphisms of $\A_n^{\bar{q},(1)}$}
\label{twisting-automorphisms}

Denote by $(1)$ the $n \times n$ skew-symmetric matrix with all entries equal to $1$.
It is straightforward to check that, for $1 \le i \le n$, there are 
$\k$-algebra automorphisms $\tau_i$ of $\A_n^{\bar{q},(1)}$ such that
\[
\begin{array}{ccrcl}
\tau_i & : & \A_n^{\bar{q},(1)} & \longrightarrow & \A_n^{\bar{q},(1)} \cr
& & x_j & \mapsto & 
\left\{ 
\begin{array}{rll}
x_j & \mbox{whenever} & 1 \le j \le i,  \cr
\lambda_{ij}^{-1}x_j & \mbox{whenever} & i < j \le n.  \cr
\end{array}
\right.
\cr
& & y_j & \mapsto & 
\left\{ 
\begin{array}{rll}
y_j & \mbox{whenever} & 1 \le j \le i, \cr
\lambda_{ij} y_j & \mbox{whenever} & i < j \le n. \cr
\end{array}
\right.
\cr
\end{array}
\]
Clearly, the above automorphisms pairwise commute. Hence, there is a morphism of groups as 
follows
\[
\begin{array}{ccrcl}
 & & \Z^n & \longrightarrow & \Aut_\k\left(\A_n^{\bar{q},(1)}\right) \cr
 & & e_i & \mapsto & \tau_i
\end{array} 
\]
It is also clear that the automorphisms $\tau_i$ are homogeneous of degree $0$ with respect to the 
$\Z^n$-grading of $\A_n^{\bar{q},(1)}$ defined in this section.

From now on, given $g\in\Z^n$, we let $\tau_g$ denote the 
$\k$-algebra automorphism -- $\Z^n$-homogeneous of degree $0$ -- 
that is the image of $g$ under the above group morphism. More precisely,
if $g=(g_1,\dots,g_n)\in\Z^n$, we have that
\[
\tau_g = \tau_1^{g_1} \dots \tau_n^{g_n}. 
\]

\subsection{Twisting of $\A_n^{\bar{q},(1)}$}\label{twisting-A}
This subsection relies on 
Zhang twists as defined in \cite{Z}. However, we are interested in left modules rather than in right ones.
So, the convenient context for us is that of \cite{RZ}; we will follow it. We will thus consider 
left twists (see \cite[Def. 1.2.1]{RZ}).

Clearly, $\tau=(\tau_g)_{g\in\Z^n}$ is a normalised left twisting system of $\A_n^{\bar{q},(1)}$ 
in the sense of \cite{RZ}.
Thus, we may associate to $\tau$ and $\A_n^{\bar{q},(1)}$ the graded $\k$-algebra  
${}^{\tau}\left(\A_n^{\bar{q},(1)}\right)$. As a graded $\k$-vector space, 
${}^{\tau}\left(\A_n^{\bar{q},(1)}\right) = \A_n^{\bar{q},(1)}$, but the algebra structure is given by a new 
associative product, that we denote $\ast$, which is defined by:
\[
y \ast z = \tau_g(y)z
\] 
whenever $z$ is an homogeneous element of degree 
$g$ of the vector space $\A_n^{\bar{q},(1)}$ and $y$ is any element of $\A_n^{\bar{q},(1)}$. 
This new product has the same unit as the 
original one and $\left({}^{\tau}\left(\A_n^{\bar{q},(1)}\right),\ast\right)$ is a $\Z^n$-graded
$\k$-algebra with respect to the original grading. 

At this stage, we have the following result.

\begin{subtheorem} --\label{TheoremA}
There is a $\k$-algebra isomorphism
$\A_n^{\bar{q},\Lambda} \longrightarrow {}^{\tau}\left(\A_n^{\bar{q},(1)}\right)$ 
such that $x_i \mapsto x_i$ and $y_i \mapsto y_i$.
\end{subtheorem}

\proof We first check that there is a $\k$-algebra homomorphism as desired. 

Let $1 \le i < j \le n$. We have that: 
\[
x_i \ast x_j - \lambda_{ij}x_j \ast x_i 
= \tau_j(x_i) x_j - \lambda_{ij} \tau_i(x_j) x_i
= x_i x_j - \lambda_{ij}\lambda_{ji}x_j x_i
= x_i x_j - x_jx_i
= 0 ;
\]
\[
y_i \ast y_j - \lambda_{ij}y_j \ast y_i 
= \tau_j^{-1}(y_i) y_j - \lambda_{ij}\tau_i^{-1}(y_j) y_i
= y_i y_j - \lambda_{ij}\lambda_{ji} y_j y_i
= y_i y_j - y_jy_i
= 0 ;
\]
\[
x_i \ast y_j - \lambda_{ji}y_j \ast x_i 
= \tau_j^{-1}(x_i) y_j - \lambda_{ji}\tau_i(y_j)x_i
= x_i y_j - \lambda_{ji}\lambda_{ij}y_j x_i
= x_i y_j - y_jx_i
= 0 ;
\]
\[
y_i \ast x_j - \lambda_{ji}x_j \ast y_i 
= \tau_j(y_i) x_j - \lambda_{ji}\tau_i^{-1}(x_j) y_i
= y_i x_j - \lambda_{ji}\lambda_{ij}x_j y_i
= y_i x_j - x_jy_i
= 0 .
\]
Now, given $1 \le i \le n$: 
\[
x_i \ast y_i - q_i y_i \ast x_i - 1 
= \tau_i^{-1}(x_i) y_i - q_i \tau_i(y_i) x_i - 1
= x_i y_i - q_i y_i x_i - 1
= 0 .
\]
This establishes the existence of such a $\k$-algebra homomorphism. 
Note that this homomorphism sends the PBW basis of the $\k$-vector space $\A_n^{\bar{q},\Lambda}$
to the PBW basis of the $\k$-vector space $\A_n^{\bar{q},(1)}$. As a consequence, it is an 
isomorphism.\qed

\subsection{Twisting of $\B_n^{\bar{q},(1)}$.}\label{twisting-B}
In this subsection, we extend 
the results of Subsection \ref{twisting-A} to the localisation $\B_n^{\bar{q},(1)}$ of 
$\A_n^{\bar{q},(1)}$.
\medskip
Let $1 \le i \le n$. In the notation of paragraph \ref{twisting-automorphisms}, it is clear that
$\tau_i(z_j)=z_j$, for all $1 \le j \le n$. It follows that the automorphism $\tau_i$ extends 
to an automorphism --still denoted $\tau_i$-- of $\B_n^{\bar{q},(1)}$. Hence, we get a 
morphism of groups 
\[
\begin{array}{ccrcl}
 & & \Z^n & \longrightarrow & \Aut_\k\left(\B_n^{\bar{q},(1)}\right) \cr
 & & e_i & \mapsto & \tau_i.
\end{array}
\] 

The same construction as in Subsection \ref{twisting-A} leads to the following result.

\begin{subtheorem} --
There is a $\k$-algebra isomorphism
$\B_n^{\bar{q},\Lambda} \longrightarrow {}^{\tau}\left(\B_n^{\bar{q},(1)}\right)$ 
such that $x_i \mapsto x_i$ and $y_i \mapsto y_i$.
\end{subtheorem}

\subsection{The $\Z^n$-grading of $P_\phi$} 

Let $\phi\in\widehat{\k[z_1^{\pm 1},\dots,z_n^{\pm 1}]}$. \\
Recall from Subsection \ref{base-de-P} that 
$P_{\phi}$ has
a $\k$-basis $\{v_k,\, k\in\Z^n\}$. Consider the $\Z^n$-grading of the $\k$-vector space
$P_\phi$ whose $k$-th component is $\k v_k$. We know, after Lemma \ref{action-xy}, that this grading
is also a grading of the $\B_n^{\bar{q},\Lambda}$-module $P_\phi$.  
Hence, we can define the twist ${}^\tau (P_\phi)$ of the 
$\B_n^{\bar{q},(1)}$-module $P_\phi$. See \cite{RZ} for the definition of
the twist of a module. It is worth mentioning at this point the following general fact:
for any $\Z^n$-graded 
$\B_n^{\bar{q},(1)}$-module $M$, $M$ and ${}^\tau M$ have the same lattice of $\Z^n$-graded 
submodules.
Using Theorem \ref{TheoremA}, we get that  
${}^\tau (P_\phi)$is naturally endowed with a structure of
$\B_n^{\bar{q},\Lambda}$-module, by restriction of scalars.\\

Now, suppose that $\bar{q}$ is generic. By Prop. \ref{generic-weight-decomposition}, the weight decomposition of $P_\phi$ is 
\[
P_\phi = \bigoplus_{\k\in\Z^n} \k.v_\k 
\]  
and each line $\k.v_k$ is the weight space of $P_\phi$ of weight $\phi\circ\sigma_k$. 
That is: the weight decomposition of $P_\phi$ and its $\Z^n$-grading coincide. 
On the other hand, any submodule of $P_\phi$ is a weight submodule 
--cf. Proposition \ref{sous-objets-et-poids}. 
It follows that any submodule or quotient of $P_\phi$ is $\Z^n$-homogeneous. As a consequence, 
the quotient $S_\phi$ is again $\Z^n$-graded and we can define its twist ${}^\tau (S_\phi)$. 
Using Theorem \ref{TheoremA}, we get that  
${}^\tau (S_\phi)$ is naturally endowed with a structure of
$\B_n^{\bar{q},\Lambda}$-module, by restriction of scalars.

\subsection{An isomorphism of $\B_n^{\bar{q},\Lambda}$-modules}\label{iso-after-twist}

In this subsection we want to consider simultaneously 
$\B_n^{\bar{q},(1)}$-modules and $\B_n^{\bar{q},\Lambda}$-modules. Hence, to avoid any ambiguity, 
we introduce the following notation. 
For all $\phi\in\widehat{\k[z_1^{\pm 1},\dots,z_n^{\pm 1}]}=R^\circ$, let
\[
P_\phi^{(1)} = \B_n^{\bar{q},(1)} \otimes_{R^\circ} \k_\phi
\qquad\mbox{and}\qquad
P_\phi^{\Lambda} = \B_n^{\bar{q},\Lambda} \otimes_{R^\circ} \k_\phi.
\]
Likewise, we denote by $\{v_\k^{(1)},\,k\in\Z^n\}$ the $\k$-basis of $P_\phi^{(1)}$ introduced in
Subsection \ref{base-de-P} and by $(v_\k^{\Lambda},\,k\in\Z^n)$ the corresponding $\k$-basis of 
$P_\phi^{\Lambda}$. Let us consider the $\k$-vector space morphism $\iota_\phi$ such that
\[
\begin{array}{ccrcl}
\iota_\phi & : & {}^\tau\left(P_\phi^{(1)}\right) & \longrightarrow & P_\phi^{\Lambda} \cr
 & & v_\k^{(1)}& \mapsto & v_\k^{\Lambda}.
\end{array}
\]
A direct application of Lemma \ref{action-xy} gives the following result.

\begin{subtheorem} --
The map $\iota_\phi$ is an isomorphism of $\B_n^{\bar{q},\Lambda}$-modules.
\end{subtheorem}

Suppose now that $\bar{q}$ is generic.
Analogously, we denote by $S_\phi^\Lambda$ the unique simple quotient 
of $P_\phi^\Lambda$ and by $S_\phi^{(1)}$ the unique simple quotient 
of $P_\phi^{(1)}$. Using the previous theorem, we have the following composition of 
$\B_n^{\bar{q},\Lambda}$-linear maps
\[
{}^\tau \left(P_\phi^{(1)}\right) 
\stackrel{\iota_\phi}{\longrightarrow} P_\phi^{\Lambda}
\stackrel{}{\longrightarrow} S_\phi^{\Lambda},
\]
where the second map is the canonical projection. This composition is surjective, hence its 
kernel is a maximal strict $\B_n^{\bar{q},\Lambda}$-submodule of   
${}^\tau \left(P_\phi^{(1)}\right)$. But the lattice of submodules of 
the ${}^{\tau}\-\-\left(\B_n^{\bar{q},(1)}\right)$-module
${}^\tau \left(P_\phi^{(1)}\right)$ coincides with 
the lattice of submodules of 
the $\B_n^{\bar{q},(1)}$-module
$P_\phi^{(1)}$. Hence, the kernel of the above map is $N_\phi^{(1)}$
--using an obvious notation.
As a consequence, we get that $S_\phi^\Lambda$ is isomorphic to 
${}^\tau \left(P_\phi^{(1)}\right)/N_\phi^{(1)}$ as a 
$\B_n^{\bar{q},\Lambda}$-module. That is, we have proved the following theorem.

\begin{subtheorem} --
There is an isomorphism of 
$\B_n^{\bar{q},\Lambda}$-modules between $S_\phi^\Lambda$ and 
${}^\tau\left(S_\phi^{(1)}\right)$.
\end{subtheorem}

We finish this section with an easy but useful observation. 

\begin{subremark} -- \rm
Let $M$ be a $\Z^n$-graded $\B_n^{\bar{q},(1)}$-module. 
Then, for all $\phi \in \widehat{\k[z_1^{\pm 1},\dots,z_n^{\pm 1}]}$, we have an equality
\[
({}^\tau M)(\phi) = M(\phi) .
\]
This is obvious since the automorphisms of the twisting system $\tau$ all act trivially on 
$z_1,\dots,z_n$.

It follows that the set of weights of ${}^\tau M$ coincides with the set of weights of $M$.
\end{subremark}

\section{Classification of simple weight modules of $\B_n^{\bar{q},\Lambda}$ in the generic case}
\label{classification}

In this final section, we collect the results of the previous ones to give an explicit description 
and a classification of the simple weight modules of the extension 
$(R^\circ,\B_n^{\bar{q},\Lambda})$ when $\bar{q}$ is generic.

We also show that the natural representation of $\B_n^{\bar{q},\Lambda}$ by skew differential 
operators on the convenient quantum affine space is a simple weight module.\\

For all this section, we fix $\bar{q}=(q_1,\dots,q_n)\in(\k^\ast)^n$ that we assume generic and a 
skew-symmetric matrix $\Lambda=(\lambda_{ij})$ with entries in $\k^\ast$. We let $\tau$ be the 
normalised left twisting system associated to $\Lambda$, as defined in Sections 
\ref{twisting-A} and \ref{twisting-B}.

As proven in Section \ref{twisting-B}, there is an algebra isomorphism 
$\B_n^{\bar{q},\Lambda} \longrightarrow {}^{\tau}\left(\B_n^{\bar{q},(1)}\right)$ 
such that $x_i \mapsto x_i$ and $y_i \mapsto y_i$. Notice that, under the above isomorphism, 
the subalgebra of 
$\B_n^{\bar{q},\Lambda}$ generated by $z_1,\dots,z_n$ and their inverses is in one-to-one correspondence with the subalgebra of ${}^\tau\left(\B_n^{\bar{q},(1)}\right)$ generated by $z_1,\dots,z_n$. Further, the latter subalgebra is in one-to-one correspondence with the subalgebra of 
$\B_n^{\bar{q},(1)}$ generated by $z_1,\dots,z_n$ by means of the identity map. 
Both subalgebras are Laurent polynomial algebras in $z_1,\dots,z_n$. By abuse of notation, we all denote them by $R^\circ$. \\

We keep the notation of Section \ref{iso-after-twist}.

\paragraph{\em Construction of simple weight modules.} 
Let $S$ be a simple weight module of $\B_n^{\bar{q},\Lambda}$.\\

By Proposition \ref{ubiquite-des-P}, there exists $\phi\in\widehat{R^\circ}$ such that
\[
S \cong S_\phi^{\Lambda}.
\]
On the other hand, we deduce from Section \ref{iso-after-twist} that 
\[
S_\phi^{\Lambda} \cong {}^\tau\left( S_\phi^{(1)}\right)
\]
and $S_\phi^{\Lambda}$ and $S_\phi^{(1)}$ have the same set of weights.

Further, by Theorem \ref{thm-simple} (see also Corollary \ref{iso-L-tenseur-L}),
if $(\phi_1,\dots,\phi_n)$ is the unique $n$-tuple of elements of $\k[z_1^{\pm 1}]$ such 
that
$\phi=\phi_1 \dots \phi_n$, then
\[
S_\phi^{(1)} \cong S_{\phi_1}^{(1)} \otimes\dots\otimes S_{\phi_n}^{(1)}.
\]
It is clear, in addition, that the set of weights of 
$S_{\phi_1}^{(1)} \otimes\dots\otimes S_{\phi_n}^{(1)}$ is the image under the map of Proposition \ref{weight-in-tpa} (2) of $w\left(S_{\phi_1}^{(1)}\right) \times\dots\times w\left(S_{\phi_n}^{(1)}\right)$, 
where, for $1 \le i \le n$, $w\left(S_{\phi_i}^{(1)}\right)$ denotes the set of weights of 
$S_{\phi_i}^{(1)}$. These last weight spaces are described in Examples 
\ref{exemple-hors-1} and \ref{exemple-1} according
to whether $\phi_i$ is or not in the $\Z$-orbit of $\1$.  \\

Denote by $\B_n^{\bar{q},\Lambda}-\wSimple$ the set of simple weight 
$\B_n^{\bar{q},\Lambda}$-modules -- note that the previous discussion shows this collection is actually a set-- and by
$\sim$ the equivalence relation on this set defined by isomorphism. We deduce from the arguments above 
that there is a surjective map
\[
\kappa_n \, : \,
\widehat{R^\circ} \longrightarrow \B_n^{\bar{q},\Lambda}-\wSimple
\longrightarrow (\B_n^{\bar{q},\Lambda}-\wSimple)/\sim
\]
that sends any $\phi\in \widehat{R^\circ}$ to the isomorphism class of the 
$\B_n^{\bar{q},\Lambda}$-module
${}^\tau \left(S_\phi^{(1)}\right)$.

\paragraph{\em Isomorphisms between simple weight modules.} Our final aim is to describe fibres of the  map $\kappa_n$.

Let $\phi,\psi\in \widehat{R^\circ}$. 
Since weight spaces and $\Z^n$-grading coincide for $S_\phi^{(1)}$ and $S_\psi^{(1)}$, 
the $\B_n^{\bar{q},(1)}$-modules $S_\phi^{(1)}$ and $S_\psi^{(1)}$ are isomorphic if and only if 
the $\B_n^{\bar{q},\Lambda}$-modules ${}^\tau\left(S_\phi^{(1)}\right)$ and 
${}^\tau\left(S_\psi^{(1)}\right)$ are isomorphic. But, the former occur if and only if
$S_\phi^{(1)}$ and $S_\psi^{(1)}$ have the same set of weights
(see, in particular, Corollary \ref{iso-S-via-poids} and Theorem \ref{thm-simple}).
It follows from this that $\phi$ and $\psi$ have the same image under $\kappa_n$ if and only if 
the corresponding simple modules $S_\phi^{(1)}$ and $S_\psi^{(1)}$ have the same set of weights.

\begin{example} -- \rm {\bf The natural representation of $\B_n^{\bar{q},\Lambda}$.} 
As already discussed, 
$\B_n^{\bar{q},\Lambda}$ may be seen as a ring of $q$-difference operators on quantum spaces. We now 
describe this representation and show it is a simple weight module. 
For this, we follow \cite[Sec. 4.1]{R2}. Let $E_n^\Lambda$ be the $\k$-algebra generated by
$y_i$, $1 \le i \le n$ subject to the relations $y_iy_j=\lambda_{ij} y_jy_i$, 
for all $1 \le i < j \le n$. This is a noetherian integral domain and we denote by $F_n^\Lambda$ its
skew-field of fractions. It is easy to check that, for $1 \le i \le n$,
there is an automorphism of $\k$-algebra as follows
\[
\begin{array}{ccrcl}
\xi_i & : & F_n^\Lambda & \longrightarrow & F_n^\Lambda \cr
 & & y_j & \mapsto & q_i^{\delta_{ij}}y_j
\end{array}
\]
and that this automorphism restrict to an automorphism of $E_n^\Lambda$.
Further, for $1 \le i \le n$, we consider the following $\k$-vector space endomorphisms of
$F_n^\Lambda$:
\[
\begin{array}{ccrcl}
m_i & : & F_n^\Lambda & \longrightarrow & F_n^\Lambda \cr
 & & f & \mapsto & y_i f
\end{array}
\quad\mbox{and}\quad
\begin{array}{ccrcl}
\partial_i & : & F_n^\Lambda & \longrightarrow & F_n^\Lambda \cr
 & & f & \mapsto & (q_i-1)y_i^{-1} (\xi_i(f)-f) 
\end{array}
\]
It is easy to see that, actually, the operators $m_i$ and $\partial_i$ stabilise $E_n^\Lambda$. 
It turns out that we have a $\k$-algebra morphism as follows:
\[
\begin{array}{ccrcl}
 & & \B_n^{\bar{q},\Lambda} & \longrightarrow & \End_\k(F_n^\Lambda) \cr
 & & x_i & \mapsto& \partial_i \cr
 & & y_i & \mapsto& m_i \cr
 & & z_i & \mapsto& \xi_i \cr 
\end{array}
\]
(see \cite[Prop. 4.1.3]{R2}) which defines an action of $\B_n^{\bar{q},\Lambda}$ on $F_n^\Lambda$
and, by restriction, an action of $\B_n^{\bar{q},\Lambda}$ on $E_n^\Lambda$. 

For $k=(k_1,\dots,k_n)\in\N^n$, let $y^k=y_1^{k_1} \dots y_n^{k_n} \in E_n^\Lambda$. As it is well
known, the family $\{y_k, k\in\N^n\}$ is a basis of the vector space $E_n^\Lambda$. Further, 
for $k\in\N^n$, it is easy to see that $y^k$ is a weight vector of weight $-k\1$, in the notation
of Section \ref{action-B} and of the introduction to Section \ref{section-exemples}. It is also
easy to verify that $E_n^\Lambda$ is simple. 

But, on the other hand, we know that the set of weights of $S_\1^{\Lambda}$ is $(-\N)^n\1$. This 
shows 
that $E_n^\Lambda$ and $S_\1^{(\Lambda)}$ are two simple weight modules of $\B_n^{\bar{q},\Lambda}$
with the same set of weights. Hence, they are isomorphic.
\end{example}


\footnotesize
\noindent V. F.:
\\Departamento de Matem\'atica\\
Universidade de S\~ao Paulo, 
S\~ao Paulo, Brasil.\\
{\tt futorny@ime.usp.br}

\medskip

\noindent L. R.:
\\Laboratoire Analyse, G\'eom\'etrie et Applications LAGA,                                                              
UMR 7539 du CNRS,\\
Institut Galil\'ee,\\
Universit\'e Paris Nord 
\\
Villetaneuse, France.\\
{\tt rigal@math.univ-paris13.fr}

\medskip
\noindent A. S.:
\\IMAS-CONICET y Departamento de Matem\'atica,
 Facultad de Ciencias Exactas y Naturales,\\
 Universidad de Buenos Aires,
\\Ciudad Universitaria, Pabell\'on 1\\
1428, Buenos Aires, Argentina. \\{\tt asolotar@dm.uba.ar}

\end{document}